\documentclass[12pt]{article}
\usepackage[utf8]{inputenc}
\usepackage{amsmath}
\usepackage{amssymb}
\usepackage{mathtools}
\usepackage{bbm}

\RequirePackage[colorlinks,citecolor=red,urlcolor=blue, linkcolor=blue]{hyperref}
\usepackage[margin=0.8in]{geometry}

\usepackage[dvips]{graphics}
\usepackage{epsfig,rotating}
\usepackage{tabularx}
\usepackage{bm}
\usepackage{bbm}

\usepackage{mdframed}
\usepackage{lipsum}
\usepackage{enumitem}

\usepackage[nottoc,numbib]{tocbibind} 
\usepackage[toc,page]{appendix}

\newtheorem{thm}{Theorem}
\newtheorem{lemma}{Lemma}

\newtheorem{definition}{Definition}

\newcommand{\beaa}{\begin{eqnarray*}}
\newcommand{\eeaa}{\end{eqnarray*}}
\newcommand{\bea}{\begin{eqnarray}}
\newcommand{\eea}{\end{eqnarray}}

\newcommand{\fvml}{\mbox{FvML}}

\newcommand{\la}{\left\{}
\newcommand{\ra}{\right\}}
\newcommand{\lb}{\left(}
\newcommand{\rb}{\right)}
\newcommand{\mb}{\mathbb}
\newcommand{\ve}{\varepsilon}
\newcommand{\argmax}{\operatornamewithlimits{argmax}}

\begin{document}
\title{Asymptotic analysis of high-dimensional uniformity tests under
heavy-tailed alternatives}
\author{Tiefeng Jiang \\ \small{School of Data Science, Chinese University of Hong Kong, Shenzhen} \\ \small{\href{mailto:jiang040@cuhk.edu.cn}{jiang040@cuhk.edu.cn}}
   \and Tuan Pham \\ \small{Department of Statistics and Data Science, University of Texas, Austin} \\  \small{\href{mailto:tuan.pham@utexas.edu}{tuan.pham@utexas.edu}}}
\date{}

\maketitle
\begin{abstract}
We study the high-dimensional uniformity testing problem, which involves testing whether the underlying distribution is the uniform distribution, given $n$ data points on the $p$-dimensional unit hypersphere. While this problem has been extensively studied in scenarios with fixed $p$, only three testing procedures are known in high-dimensional settings: the Rayleigh test \cite{Cutting-P-V}, the Bingham test \cite{Cutting-P-V2}, and the packing test \cite{Jiang13}. Most existing research focuses on the former two tests, and the consistency of the packing test remains open. We show that under certain classes of alternatives involving projections of heavy-tailed distributions, the Rayleigh test is asymptotically blind, and the Bingham test has asymptotic power equivalent to random guessing. In contrast, we show theoretically that the packing test is powerful against such alternatives, and empirically that its size suffers from severe distortion due to the slow convergence nature of extreme-value statistics. By exploiting the asymptotic independence of these three tests, we then propose a new test based on Fisher's combination technique that combines their strengths. The new test is shown to enjoy all the optimality properties of each individual test, and unlike the packing test, it maintains excellent type-I error control. As a consequence of the asymptotic independence, we derive the non-null limiting distribution of the packing test without using the Chen-Stein method for Poisson approximation.
\end{abstract}
 \tableofcontents

\section{Introduction}

Uniformity testing stands as one of the most crucial problems in modern directional statistics. Consider the hypersphere $\mb{S}^{p-1}= \la x \in \mb{R}^{p}: \|x\|_2 =1 \ra$, where $\|.\|_2$ is the Euclidean distance. The observed data points are denoted by $\bm{X}_1, \bm{X}_2,...,\bm{X}_n$ with $\bm{X}_i \in \mb{S}^{p-1}$ for all $i=1,2,...,n$.  We assume that $\bm{X}_i$'s are drawn independently from an unknown distribution $\mu$ supported on the hypersphere $\mb{S}^{p-1}$.  The uniform distribution on the hypersphere $\mb{S}^{p-1}$ is denoted by $\mbox{Unif}(\mb{S}^{p-1})$. The uniformity testing problem can be formulated as
\begin{align} \label{uniform-test}
    H_0: \mu= \mbox{Unif}(\mb{S}^{p-1}) \ \ \ \mbox{against}\ \ \ H_{1}:\mu \neq \mbox{Unif}(\mb{S}^{p-1}).
\end{align} 
In fixed dimensions, particularly when $p=2$ and $p=3$, uniformity testing holds significant importance across various applications in fields such as geology, paleomagnetism, and cosmology. For a detailed treatment of these topics and their applications, readers are encouraged to refer to the monographs \cite{Fisher, Ley-Verdebout, M-Jupp}. Recent tests that work for arbitrary but fixed dimensions include \cite{garcia2023projection,fernandez2023new,garcia2021cramer} (see also the references therein). The readers are referred to the survey papers \cite{survey-uni,pewsey2021recent} for recent progress on this problem and for a list of recent testing procedures.

In the era of big data, there has been growing interest in studying high-dimensional directional statistics. One common application of uniformity tests is in testing for spherical and elliptical symmetry, a fundamental question in multivariate statistical modeling. A random vector \(\bm{X} \in \mathbb{R}^p\) is spherically symmetric if and only if \(\|\bm{X}\|\) is independent of \(\bm{X}/\|\bm{X}\|\), and \(\bm{X}/\|\bm{X}\|\) is uniformly distributed on \(\text{Uni} \left( \mathbb{S}^{p-1} \right)\). Thus, spherical symmetry can be tested by combining independence tests and uniformity tests. In practice, violations of spherical symmetry often result from a failure to fulfill uniformity on the hypersphere for projected data (see \cite{Cutting-P-V} for details). Recent developments in this area can be found in \cite{albisetti2020testing, tang2023nonparametric}, along with the references therein.

Sign-based procedures are also widely used in shape analysis and nonparametric statistics. This approach projects observations onto hyperspheres and conducts statistical inference based on the projected data. Its robustness in high dimensions is due to the concentration of measure phenomenon, which suggests that the directions, rather than the magnitudes, of observations carry the majority of information. Below, we provide a brief, non-exhaustive overview of the high-dimensional directional statistics literature.

In \cite{Dryden05}, the author explores the asymptotic properties of high-dimensional spherical distributions and their applications in brain shape modeling. Clustering analysis on large-dimensional hyperspheres has been studied in \cite{Banerjee04, Banerjee03}. The potential applications of high-dimensional uniformity tests were demonstrated in \cite{Juan2001}, where the authors relate the problem of multivariate outlier detection to uniformity testing. Sign-based procedures in high dimensions have been investigated in \cite{Zou14} in the context of sphericity testing with unspecified locations and in \cite{WPL15}, where the authors propose a high-dimensional nonparametric mean test. Problems involving testing for concentration parameters in high-dimensional settings have been examined in \cite{Ley-P-2}.

Despite the extensive literature on fixed-dimensional tests, the uniformity testing problem in the high-dimensional context with diverging dimensions remains relatively unexplored. As the dimension diverges to infinity, many existing procedures require significant adjustments to function properly. Additionally, there is often no tractable limiting distribution under uniformity, and the power of these tests tends to be low due to the curse of dimensionality. To the best of the authors' knowledge, only three high-dimensional tests have been investigated in the literature. A brief overview of these tests is provided below.
\begin{enumerate}
    \item {\it Rayleigh test in \cite{Cutting-P-V} and \cite{Ley-P}}. This test can be formulated in terms of a U-statistic of the data points with the inner product kernel, i.e.
     \begin{align}
             R_n &:= \frac{\sqrt{2p}}{n} \sum_{1 \leq i<j \leq n} \bm{X}^{\top}_i \bm{X}_j. \label{Rayleigh}
     \end{align}
    \item {\it Bingham test in \cite{Cutting-P-V2, Zou14} and \cite{Ley-P}}. This test is also based on a U-statistic of the data points, but with a quadratic inner product kernel, i.e. 
    \begin{align}
            B_n &:= \frac{p}{n} \sum_{1 \leq i<j \leq n} \Big[ \lb \bm{X}^{\top}_i \bm{X}_j \rb^2 - \frac{1}{p}  \Big]. \label{Bingham}
    \end{align}

    \item {\it Packing test in \cite{Jiang13}.} This test is based on the smallest angle, i.e.
    \begin{align}
            P_n &:= p \cdot  \max_{1 \leq i<j \leq n} \lb \bm{X}^{\top}_i \bm{X}_j \rb^2 -4 \log n + \log \log n \label{packing}.
    \end{align}
\end{enumerate}
It is known that the Rayleigh test $R_n$ and the Bingham test $B_n$ enjoy a \textit{doubly robust} property: under the null hypothesis and the single assumption $\min \la n,p \ra \to \infty$, both $R_n$ and $B_n$ converge in distribution to the standard normal distribution. This feature is highly desirable since no restriction on the dependence between $p$ and $n$ is imposed, and neither resampling procedures nor tuning parameters are needed to get the critical values of such tests. Regarding the packing test $P_n$, it is known that, under the null hypothesis and the mild assumption $p \gg (\log n)^2$, $P_n$ converges in distribution to the Gumbel distribution with CDF $\exp \lb -(8\pi)^{-1/2} e^{-x/2} \rb$; see \cite{Jiang13}. 

The power analysis of these three tests in high dimensions, however, remains extremely limited. It has been shown in \cite{Cutting-P-V} that the Rayleigh test is asymptotically optimal against the class of Fisher-von Mises-Langevin (FvML) distributions, and in \cite{Cutting-P-V2} that the Bingham test is asymptotically optimal against the class of Watson distributions. Existing results for both the Rayleigh and Bingham tests primarily focus on local alternatives, with the main tool being Le Cam's theory and martingale central limit theorem.  

Regarding the packing test, there has been no result concerning its consistency properties, although simulation studies suggest that its power is low against the aforementioned classes of alternatives. The analysis of the packing test is particularly challenging, as the martingale central limit theorem does not apply to extreme-value statistics. Typically, the optimal tests derived from LAN expansions are often U-statistics, whose asymptotic distribution is Gaussian under local alternatives and can be analyzed via the martingale central limit theorem. However, extreme-value statistics, such as the packing test, are highly nonlinear, and their asymptotic distributions are unstable under small perturbations in the dependent structures. For example, the largest entries of high-dimensional sample correlation matrices with an autoregressive correlation structure driven by a vanishing sequence of parameters \(\rho_n\) might not asymptotically follow the Gumbel law, whereas the Gumbel law can still hold in the case of sphericity (\(\rho_n \equiv 0\)) or when \(\rho\) is a large constant (\(\rho_n \equiv c \in (0,1)\)); see \cite{jiang2023asymptotic} for details.

Since it is well understood that the Rayleigh and Bingham tests possess strong optimal properties for testing against the class of FvML and Watson distributions, with explicit asymptotic power within these classes (see \cite{Cutting-P-V, Cutting-P-V2}), we shall not attempt to analyze their power under local alternatives in this paper. However, a natural question that arises from previous analyses is: in what scenarios does the packing test outperform the other two tests? Furthermore, how can we combine their strengths to develop a new procedure that is more robust and powerful against a boarder range of alternatives?

The key finding of this article is that the packing test performs exceptionally well against heavy-tailed alternatives, whereas the other two tests do not. This demonstrates that there exist alternatives that are "far apart" from the uniform distribution but, due to the curse of dimensionality, their realizations spread out similarly to the uniform distribution and do not exhibit any noticeable clusters. As Theorem \ref{Packing-asymp} suggests, random points drawn from heavy-tailed alternatives behave differently from their uniform counterparts: while most points remain far apart, a few are either very close to one another or form nearly straight lines. Heuristically, the Rayleigh and Bingham tests can only detect one or two large clusters in the data, which explains why they fail to be consistent—since the majority of random points under \(\mu_{\alpha,p}\) do not form any clusters. Such alternatives are difficult to detect using classical tests, which is why we propose combining the strengths of multiple procedures to develop novel tests that are powerful against a broader range of alternatives. Specifically, our contributions are:

    

{\indent (i)} We analyze the asymptotic distribution of $R_n$ and $B_n$ under projections of heavy-tailed distributions (see Definition \ref{alpha-spherical} below)  and show that $R_n$ is asymptotically blind to such alternatives, while $B_n$ has an asymptotic power of $0.5$, which is equivalent to random guessing. 

{\indent (ii)}  In contrast to the two tests $R_n$ and $B_n$, we show that the packing test $P_n$ is consistent under the same class of alternatives. In particular, we prove a law of large numbers for the largest off-diagonal entries of sample correlation matrices with heavy-tailed entries. This result can be of independent interest. 

{\indent (iii)} Although the test $P_n$ is powerful against heavy-tailed alternatives, simulation studies suggest that its size is severely distorted. By showing that $R_n$, $B_n$, and $P_n$ are asymptotically independent, we propose a new test based on Fisher's combination technique that enjoys the best of all worlds: it is powerful against all known alternatives while having its empirical size close to the prescribed nominal level.

{\indent (iv)} As a byproduct of the asymptotic independence result, we illustrate how to derive the non-null asymptotic distribution of the packing test via Le Cam's third lemma. This approach, being of independent interest, is able to bypass the need for the Chen-Stein method for Poisson approximation, which is known to be complicated when the dependence is complex. We believe this technique is applicable to other extreme-value based test statistics and also in other settings.

The rest of the article is organized as follows. In Section \ref{sec-main}, we introduce the class of heavy-tailed alternatives considered in this paper and state our main results. We present the new Fisher's combination test and discuss the consequence of the asymptotic independence in Section \ref{sec-power enhancement}. We perform simulation studies in Section \ref{sec-sim}. Some conclusions, remarks, and directions for future research are discussed in Section \ref{sec-discuss}. The proofs of the main results are presented in Section \ref{sec-proofs} and some technical lemmas can be found in Section \ref{sec-technical}.

\section{Models and main results} \label{sec-main}

\subsection{Heavy-tailed alternatives}

 In this section, we use the standard terminologies of heavy-tailed modelling, which can be found, for example, in the recent monograph \cite{Kulik2020}. Let $X_1, X_2, \cdots, X_p$ be i.i.d. random variables with regularly varying tail of index $\alpha \in (0,2)$. The last condition means that
\begin{align} \label{reg-varying}
   \mb{P} \lb |X_1| \geq x \rb \sim L(x) x^{-\alpha}
\end{align}
as $x \to \infty$, for some function $L(x)$ such that
$$\lim_{x \to \infty} \frac{L(ax)}{L(x)} =1$$
for all $a \in \mb{R}$. Such functions $L$ are called ``slowly varying function".
It is also well-known that the distributions of the form \eqref{reg-varying} have infinite variance for all $\alpha \in (0,2)$ and infinite first moment for $\alpha \in (0,1)$. Heavy-tailed alternatives are ubiquitous in high-dimensional datasets. For example, it is believed that the underlying distributions of gene expression levels are heterogeneous, heavy-tailed, and exhibit complex dependence structures (see \cite{CH09}). In fact, heavy-tailed alternatives have been considered in the context of sphericity testing, a problem closely related to \eqref{uniform-test} (see, for example, \cite{Zou14} and the references therein).

The class of alternatives we will be dealing with is the projections of $\lb X_1, X_2,\cdots, X_p\rb$ onto the hypersphere $\mb{S}^{p-1}$, which we shall refer to as {\it $\alpha$-spherical distributions}. The precise definition can be formulated as follows

\begin{definition}[$\alpha$-spherical distributions] \label{alpha-spherical}
We say that a distribution $\mu_{\alpha, p}$ supported on $\mb{S}^{p-1}$ is  an $\alpha$-spherical distribution if it satisfies
\begin{align*} 
    \mu_{\alpha,p} \stackrel{d}{=} \lb \frac{X_1}{\sqrt{\sum_{k=1}^{p} X_k^2}}, \frac{X_2}{\sqrt{\sum_{k=1}^{p} X_k^2}},\cdots, \frac{X_p}{\sqrt{\sum_{k=1}^{p} X_k^2}} \rb
\end{align*}
where $X_i$'s are i.i.d. with common law $\mu$ and $\mu$ is regularly varying with index $\alpha \in (0,2)$ in the sense of \eqref{reg-varying}. Additionally, if $\mu$ is symmetric, we say that $\mu_{\alpha,p}$ is a symmetric $\alpha$-spherical distribution.
\end{definition}
Although the $X_i$'s in Definition \ref{alpha-spherical} above have infinite variance, the coordinates of $\mu_{\alpha}$ are bounded and possess finite moments of arbitrary order. If one makes the extra assumption that the corresponding law $\mu$ of $X_1$ is symmetric, then all the coordinates of $\mu_{\alpha,p}$ are centered and have variance $1/p$, which is the same as that of $\mbox{Unif}\lb \mb{S}^{p-1} \rb$. It is therefore not surprising that, when $p$ is large, $\mu_{\alpha,p}$ looks very similar to the uniform distribution and does not display any sort of axial or non-axial  patterns. This heuristic explains why both the Rayleigh test and the Bingham test fail to be consistent against such alternatives since they can only detect either axial or non-axial patterns in the data. 

The key difference between \(\mu_{\alpha,p}\) and the uniform distribution lies in their mixed moments of orders higher than four, as reflected in the formula \eqref{heavy-mixed-monent}. Specifically, this implies that 
$$
\mathbb{E} \left[ X_{i_1}^{2k_1} X_{i_2}^{2k_2} \cdots X_{i_r}^{2k_r} \right] \sim \frac{\text{constant}(\alpha)}{p^r}
$$
as \(p \to \infty\), where \((X_1, X_2, \dots, X_p)\) is drawn from \(\mu_{\alpha,p}\). In contrast, for the uniform distribution, one has
$$
\mathbb{E} \left[ X_{i_1}^{2k_1} X_{i_2}^{2k_2} \cdots X_{i_r}^{2k_r} \Big| \text{Uni} \left( \mathbb{S}^{p-1} \right) \right] \sim \frac{\text{constant}}{p^{\sum_{i=1}^{r} k_i}}
$$
as \(p \to \infty\). Thus, the mixed moments of \(\mu_{\alpha,p}\) decay much more slowly than those of \(\text{Uni} \left( \mathbb{S}^{p-1} \right)\). This observation has been used to derive the limiting empirical spectral distribution (ESD) of sample correlation matrices with heavy-tailed distributions in \cite{heiny2022}, where the authors proved that the limiting ESD converges to a new type of distribution that continuously extends the Marchenko-Pastur law. This fact will also play an important role in the proofs of Theorems \ref{Rayleigh-asymp} and \ref{Bingham-asymp}.

Additionally, we can prove that, with high probability, there exist two data points such that the angle between them is either close to \(0\) or close to \(\pi\). This behavior is fundamentally different from the uniform distribution, where in high dimensions, all pairs of points are almost always orthogonal (see \cite{Jiang13} for details). As a result, the

\subsection{Asymptotic power of the Rayleigh and Bingham tests}
Let us start with the result concerning the asymptotic distribution of the Rayleigh test under $\mu_{\alpha,p}$. Recall the definition of  $\alpha$-spherical distributions from Definition \ref{alpha-spherical}. In what follows, we assume that the data $\bm{X}_1$, $\bm{X}_2,\cdots,\bm{X}_n$ are drawn independently from a symmetric $\alpha$-spherical distribution $\mu_{\alpha,p}$ with $\alpha \in (0,2)$. The symmetry assumption is imposed in Theorem \ref{Rayleigh-asymp} and \ref{Bingham-asymp} to get exact expressions of the asymptotic distribution. In the absence of the symmetric assumption, it is likely that the asymptotic distribution is still normal, but with different scaling. However, we do not pursue this technical consideration here. 

It is worth noting that our results in Theorem \ref{Rayleigh-asymp} and \ref{Bingham-asymp} below do not follow from the results proven in \cite{Cutting-P-V} and \cite{Cutting-P-V2} since the $\alpha$-spherical alternatives we consider here do not belong to the classes of axial or non-axial distributions considered in their papers. Let us start with the asymptotic power of the Rayleigh test.

\begin{thm} \label{Rayleigh-asymp}
    Let $p=p_n$ such that $p \to \infty$ and $p/n^2 \to 0$. Recall $R_n$ from \eqref{Rayleigh}. Then, under $\mu_{\alpha, p_n}$, $R_n$ converges in distribution to a standard normal distribution. 
\end{thm}
An immediate consequence of Theorem \ref{Rayleigh-asymp} is that the Rayleigh test is blind to the symmetric $\alpha$-spherical alternatives $\mu_{\alpha,p_n}$, for any $\alpha \in (0,2)$. This is because the asymptotic distribution of $R_n$ is the same under both uniformity and $\mu_{\alpha,p_n}$. Therefore, the Rayleigh test with size $\beta$ also has power asymptotically equal to $\beta$. This agrees remarkably well with the simulation studies in Section \ref{sec-sim}. Next, we consider the Bingham test.

\begin{thm} \label{Bingham-asymp}
    Suppose $p/n \to \gamma \in (0,\infty)$. Recall $B_n$ in \eqref{Bingham}. Then, under $\mu_{\alpha, p_n}$, we have
    $$\frac{\sqrt{n}}{p} B_n \stackrel{d}{\to} N \lb 0, \frac{(2 - \alpha)^2}{8 \gamma} \rb.$$
\end{thm}
From Theorem \ref{Bingham-asymp}, one can compute the asymptotic power against $\mu_{\alpha, p_n}$ as follows. Suppose the Bingham test at nomimal level $\beta$ rejects if $B_n \geq q_{\beta}$, where $q_{\beta}$ is the $(1-\beta)$-quantile of the standard normal distribution. Then, the asymptotic power under $\mu_{\alpha,p_n}$ is 
\begin{align*}
\lim_{n \to \infty} \mb{P}_{\mu_{\alpha, p_n}} \lb B_n \geq q_{\beta} \rb &= \lim_{n \to \infty} \mb{P}_{\mu_{\alpha, p_n}} \lb \frac{\sqrt{n}}{p} B_n \geq q_{\beta} \cdot \frac{\sqrt{n}}{p} \rb \\
&=  \mb{P} \lb N\lb 0,\frac{(2 - \alpha)^2}{8 \gamma} \rb \geq 0 \rb = 0.5. 
\end{align*}
Thus, under $\mu_{\alpha,p_n}$ and the regime $p/n \to \gamma$, $B_n$ has asymptotic power equivalent to random guessing, regardless of the nominal size $\beta$.  The simulation studies in Section \ref{sec-sim} also confirm our theoretical findings. 

\subsection{Consistency of the packing test}

In this subsection, we assume that the data $\bm{X}_1$, $\bm{X}_2,\cdots,\bm{X}_n$ are drawn independently from a $\alpha$-spherical distribution $\mu_{\alpha,p}$ with $\alpha \in (0,2)$ (see Definition \ref{alpha-spherical} above). We do not assume symmetry in what follows. Our main result is
\begin{thm} \label{Packing-asymp}
    Suppose $p \to \infty$ such that $p=o(n^2)$. Recall $P_n$ from \eqref{packing}. Then, under $\mu_{\alpha,p}$, we have 
    $$\max_{1 \leq i<j \leq n} \lb \bm{X}^{\top}_i \bm{X}_j \rb \stackrel{\mb{P}}{\to} 1$$
    as $n \to \infty$. Here $\stackrel{\mb{P}}{\to}$ indicates the convergence in probability.
\end{thm}
Theorem \ref{Packing-asymp} immediately yields the consistency of the packing test since under uniformity, we have (see \cite{Jiang13}, Theorem 6)
$$\max_{1 \leq i<j \leq n} \lb \bm{X}^{\top}_i \bm{X}_j \rb^2 =  \frac{4 \log n}{p} + o_{\mb{P}}(1).$$
o the best of our knowledge, this is the first result addressing the consistency of the packing test. Simulation results also suggest that the packing test performs poorly against common alternatives, such as the FvML distributions or the Watson distributions, though no rigorous proofs have been provided so far. Our result in Section \ref{consequence} also provides justification for why the packing test does not have non-trivial power against FvML alternatives.

Interestingly, the largest inner product in Theorem \ref{Packing-asymp} is the same as the largest off-diagonal entries of the sample correlation matrices with regularly-varying distributions of index $\alpha \in (0,2)$. Equivalently, Theorem \ref{Packing-asymp} asserts that the largest off-diagonal entries of sample correlation matrices converge to 1 in probability under heavy-tailed alternatives. This is also a new result, which differs significantly from the light-tailed settings. It is worth noting that sample correlation matrices under heavy-tailed distributions have remained unexplored in the literature until the recent line of works \cite{heiny2017eigenvalues,heiny2019eigenstructure,heiny2021large,heiny2022}, in which the authors show that sample covariance and correlation matrices behave much differently under heavy-tailed distributions than their light-tailed counterparts.

Theorem \ref{Packing-asymp} also provides insights into the geometry of the set of random inner products $\la \bm{X}^{\top}_i \bm{X}_j; 1 \leq i < j \leq p \ra$. It is straightforward to verify that \(\mathbb{E} \left[\bm{X}^{\top}_i \bm{X}_j\right] = 0\) and \(\text{Var}\left[\bm{X}^{\top}_i \bm{X}_j\right] = 1/p\) for any \(i \neq j\). Therefore, on average, most of the points are orthogonal to each other when $p$ is large. However, with high probability, there will always be two points that are either very close to each other or form an almost straight line. This behavior is entirely different from that of random points sampled under the uniform distribution, as studied in \cite{Jiang13}, where the points are almost always pairwise orthogonal.

The proof of Theorem \ref{Packing-asymp} does not follow from the typical self-normalized large deviation theory commonly used in the light-tailed case. In fact, such a result is not available in the literature for heavy-tailed distributions. The only relevant result in the heavy-tailed case that we are aware of is \cite{shao1997self}, which provides large and moderate deviation approximations to the tail probability
$$
\mathbb{P} \left( \frac{\sum_{k=1}^{p} X_k}{\sqrt{\sum_{k=1}^{p} X_k^2}} \geq x \right)
$$
as \(x \to \infty\), where the \(X_i\)'s are i.i.d. regularly varying random variables with index \(\alpha \in (0,2)\). However, this result cannot be adapted to our setting because the nature of the problem is fundamentally different: the fraction in the probability above has a sub-Gaussian limiting distribution as \(p \to \infty\), while the inner product under \(\mu_{\alpha,p}\) behaves asymptotically like the ratio of two independent \(\alpha\)-stable distributions, which also exhibit heavy tails.


\section{Power enhancement} \label{sec-power enhancement}

\subsection{Fisher's combination technique}

In the previous section, we have seen that the packing test outperforms the other two tests when testing against heavy-tailed alternatives. Unfortunately, the packing test is not recommended in practice due to its severe size distortion, which is a consequence of extreme-value statistic's slow convergence. We therefore propose a new test that inherits the strength of all three tests while maintaining excellent empirical sizes. Our main result in this section is

\begin{thm} \label{asymp-inde}
    Let $R_n, B_n$ and $P_n$ be defined in (\ref{Rayleigh}), (\ref{Bingham}) and (\ref{packing}), respectively. Then, under uniformity,  for any fixed $(x,y,z) \in \mb{R}^3$, we have 
    $$ \mb{P} \lb R_n \leq x, B_n \leq y, P_n \leq z \rb \to \Phi(x) \cdot \Phi(y) \cdot G(z),$$
    provided $ p/(\log n)^2 \to \infty$. Here $\Phi(x)$ is the CDF of a standard normal distribution and $G(x) = \exp \lb -(8\pi)^{-1/2} e^{-x/2} \rb$ is the limiting distribution of $P_n$.
\end{thm}
Theorem \ref{asymp-inde} asserts that the three tests $R_n, B_n$ and $P_n$ are asymptotically independent. Based on this, we can define a new $\beta$-level test, called the Fisher's combination test $\phi^F_n$ as
\begin{align} \label{Fisher-test}
  \phi^{F}_n =  \mathbf{1}_{\la C_n \leq 1-(1-\beta)^{1/3} \ra}
\end{align}
where 
    $$C_n =  \min \la 1 - \Phi(R_n), 1- \Phi(B_n), 1- G(P_n) \ra.$$ 
In other words, the test $\phi^F_n$ rejects if the smallest $p$-value of the three tests $R_n, B_n$ and $P_n$ is below a certain threshold. The Fisher's combination techniques in high-dimensional statistics were initially introduced in a series of works by \cite{ Yu1, Yu2, Yu3} for testing mean and covariance structures. It was also employed in the context of testing for cross-sectional dependence in panel models by \cite{Feng}. However, to the best of our knowledge, it has not been used in the context of directional data before. The proof of Theorem \ref{asymp-inde} is based on a leave-one-out type argument combined with the inclusion-exclusion principle.

Note that the test $\phi^F_n$ asymptotically achieves the desired size $\beta$ in the sense that
$$\lim_{n \to \infty} \mb{P} \lb  \phi^{F}_n = 1 \Big| \mbox{Unif} \lb \mb{S}^{p-1} \rb\rb = \beta.$$
Therefore, it is much less conservative than the popular Bonferroni procedure. It is also clear that the proposed test $\phi^F_n$ is rate-optimal or consistent against any class of alternatives if at least one of the three tests $R_n$, $B_n$ or $P_n$ is. Based on the results proved in \cite{Cutting-P-V} and \cite{Cutting-P-V2}, we can deduce that the test in \eqref{Fisher-test} is minimax-optimal against the class of FvML distributions and the class of Watson distributions. The optimality also holds for the semiparametric models considered in \cite{Cutting-P-V,Cutting-P-V2} that includes the FvML distributions and the Watson distributions as special cases. Moreover, thanks to Theorem \ref{Packing-asymp}, it is robust to heavy-tailed alternatives, which is a type of alternatives that neither the Rayleigh test nor the Bingham test performs well.  We demonstrate the desirable empirical performance of the test $\phi^F_n$ in the simulation studies conducted in Section \ref{sec-sim} below.

\subsection{Consequence of the asymptotic independence} \label{consequence}

In this section, we illustrate an interesting consequence of Theorem \ref{asymp-inde}: via Le Cam's third lemma, we derive the asymptotic distribution of the packing test under the family of FvML distributions. Unlike U-statistics-based tests, such as the Rayleigh and Bingham tests, the packing test is of the extreme-value type and lacks martingale structures. This absence makes it significantly more challenging to analyze the non-null asymptotic distributions of the packing test. Existing approaches often depend on the Chen-Stein method for Poisson approximation, but this technique becomes difficult to apply when the dependence among the random variables contributing to the maximum is complex.

The argument outlined below is based on the asymptotic independence result in Theorem \ref{asymp-inde} and does not rely on the Chen-Stein method. We believe it can be generalized in a more systematic way, but additional machinery is needed to achieve such a framework. Nevertheless, we demonstrate this with an example involving the FvML distributions. Recall that the class of FvML distributions is given by
\begin{align} \label{fvml}
\frac{d\text{FvML}(\kappa, \bm{\mu})}{d \mu_0 } \lb \bm{x} \rb \propto \exp \lb \kappa \bm{\mu}^{\top} \bm{x} \rb,
\end{align}
where $\mu_0:=\text{Uni} \lb \mathbb{S}^{p-1} \rb$, $\kappa > 0$ is the concentration parameter and $\bm{\mu}$ is the location parameter.

Consider the regime where $\kappa_n = \tau_n p^{3/4}/\sqrt{n}$ and $\tau_n \to \tau>0$. It is known from \cite{Cutting-P-V} that the rate $p^{3/4}/\sqrt{n}$ is optimal for the problem \eqref{uniform-test} within the class of FvML distributions, and this rate is achieved by the Rayleigh test $R_n$. Let $\mathcal{T}_n$ be the maximal invariant of the orthogonal group $\mbox{SO}(p_n)$ consisting of $p_n \times p_n$ orthogonal matrices. Theorem 4.2 from \cite{Cutting-P-V} then gives the LAN expansion
\begin{align} \label{LAN}
\log \frac{d \mb{P}^{(\mathcal{T}_n)}_{\mu_n}}{d \mb{P}_{\mu_0}} =  \frac{\tau_n^2}{\sqrt{2}} R_n - \frac{\tau_n^4}{4} + R^{'}_n 
\end{align}
where $\mu_n:= \mbox{FvML}(\kappa_n, \bm{\mu}_n)$, $R_n$ is the Rayleigh test as defined in \eqref{Rayleigh} and $R^{'}_n = o_{\mb{P}}(1) $.
Thanks to Theorem \ref{asymp-inde} and the assumption $\tau_n \to \tau$, we get the joint limiting distribution
\begin{align*} 
    \lb P_n,  \log \frac{d \mb{P}^{(\mathcal{T}_n)}_{\mu_n}}{d \mb{P}_{\mu_0}}\rb 
  &=  \lb P_n, \frac{\tau_n^2}{\sqrt{2}} R_n - \frac{\tau_n^4}{4} \rb + \lb 0,  R^{'}_n\rb \\
  \stackrel{d}{\to} & \lb G, N \lb -\frac{\tau^4}{4} , \frac{\tau^4}{2} \rb \rb.
\end{align*}
Here we use the fact that $R^{'}_n = o_{\mb{P}}(1) $ to deduce the last line with $G$ being the Gumbel distribution appears in the asymptotic distribution of $P_n$. Note that $G$ and the normal distribution in the display above are independent due to Theorem \ref{asymp-inde}. Apply the continuous mapping lemma with the function $(x,y) \mapsto (x,e^y)$, we arrive at
\begin{align*}
    \lb P_n , \frac{d \mb{P}^{(\mathcal{T}_n)}_{\mu_n}}{d \mb{P}_{\mu_0}} \rb \stackrel{d}{\to} \lb G, \exp \lb N \lb -\frac{\tau^4}{4} , \frac{\tau^4}{2} \rb \rb \rb.
\end{align*}
 By using  Le Cam's third lemma,  the convergence above implies that, under the sequence of alternatives $\fvml( \kappa_n, \bm{\mu}_n)$, $P_n$ converges in distribution to the Gumbel law with CDF $G$. Here the maximal invariant $\mathcal{T}_n$ does not affect the distribution since $P_n$ is rotationally invariant.

 In other words, $P_n$ has the same asymptotic distribution under $\fvml( \kappa_n, \bm{\mu}_n)$ as it does under $H_0$. This also implies that the packing test is blind to such alternatives at the detectable threshold $p^{3/4}/\sqrt{n}$. The heart of the argument above is the LAN expansion \eqref{LAN}.  Whenever such an expansion exists, it is likely that we can apply a similar argument by establishing an appropriate result concerning the asymptotic joint distribution of the packing test and the statistic that appears in the LAN expansion, which is often a U-statistic with a tractable asymptotic distribution.

\section{Simulation studies} \label{sec-sim}

We perform a Monte Carlo simulation study to verify the size and power of our proposed testing procedure. Let us start with the size's simulation. Recall the tests  $\phi^F_n$, $R_n$, $B_n$ and $P_n$ defined in \eqref{Fisher-test}, \eqref{Rayleigh}, \eqref{Bingham} and \eqref{packing}, respectively. We will compare their sizes of the tests  in three scenarios: $(n,p)=(80,100)$, $(100,100)$ and $(100,120)$. The results are summarized in Table \ref{Table-size} below.

\begin{table}[ht]
\def~{\hphantom{0}}
\caption{Empirical size at $\alpha=0.05$}{%
\begin{tabular}{ |c|c|c|c| }
\hline
\textsf{Test/Dimension} & \textsf{n=80, p=40} & \textsf{n=100,p=100} & \textsf{n=100,p=120} \\
\hline
$\phi^F_n$ & 0.06  & 0.045  & 0.047 \\
\hline
$R_n$ & 0.075   &  0.0495  & 0.0515\\ 
\hline
$P_n$ & 0.132 & 0.111 & 0.103\\
\hline
$B_n$ & 0.0477 & 0.0485  & 0.045 \\
\hline
\end{tabular}}
\label{Table-size}
\end{table}

Based on Table \ref{Table-size}, it is evident that the proposed test $\phi^F_n$, Rayleigh test $R_n$, and the Bingham test $B_n$ demonstrate effective control over the type I error rates, particularly for values of $p$ equal to or greater than $100$. The Bingham test $B_n$ exhibits robust performance even for small sample sizes and dimensions. As expected, the packing test $P_n$ performs poorly under the selected sample sizes and dimensions, attributed to the slow convergence rate inherent in extreme value statistics. 

For the power simulation, we set $\mu$ (see Definition \ref{alpha-spherical}) to have the following three types of distributions: the chi-squared distribution with one degree of freedom (centered and normalized to have variance one), the Cauchy distribution and the $t$-distribution with $1.5$ degree of freedom. Note that in the first setting, the distribution of the coordinates is of exponential tail and is not covered by our theories while the other two settings do. The results are presented in Table \ref{Table-power} below. 

\begin{table}[ht]
\def~{\hphantom{0}}
\caption{Empirical power under heavy-tailed alternatives}{%
\begin{tabular}{ |c|c|c|c|c| }
\hline
\textsf{Test/Dimension} & \textsf{Distribution of $\mu$} & \textsf{n=80, p=40} & \textsf{n=100,p=100} & \textsf{n=100,p=120} \\
\hline
$\phi^F_n$ & \begin{tabular}{ccc}
  $ \chi^2(1)$   \\
   $ \mbox{Cauchy}$ \\
   $t_{1.5}$
\end{tabular} & \begin{tabular}{ccc}
    0.357  \\
     1 \\
     0.9995
\end{tabular}  & \begin{tabular}{ccc}
   0.6125   \\
    1 \\
    1
\end{tabular}  & \begin{tabular}{ccc}
     0.6245 \\
     1 \\
     1
\end{tabular} \\ 
\hline
$R_n$  & \begin{tabular}{ccc}
  $ \chi^2(1)$   \\
   $ \mbox{Cauchy}$ \\
   $t_{1.5}$
\end{tabular} & \begin{tabular}{ccc}
     0.134  \\
     0.056 \\
     0.063
\end{tabular}  &  \begin{tabular}{ccc}
    0.0695 \\
    0.0515 \\
    0.0555
\end{tabular}  & \begin{tabular}{ccc}
    0.0755  \\
     0.059 \\
     0.057
\end{tabular}\\ 
\hline
$B_n$  & \begin{tabular}{ccc}
  $ \chi^2(1)$   \\
   $ \mbox{Cauchy}$ \\
   $t_{1.5}$
\end{tabular} & \begin{tabular}{ccc}
     0.104  \\
     0.225 \\
     0.1745
\end{tabular}  &  \begin{tabular}{ccc}
    0.073 \\
    0.3105 \\
    0.242
\end{tabular}  & \begin{tabular}{ccc}
     0.0715  \\
    0.339 \\
    0.27
\end{tabular}\\
\hline 
$P_n$  & \begin{tabular}{ccc}
  $ \chi^2(1)$   \\
   $\mbox{Cauchy}$\\
   $t_{1.5}$
\end{tabular} & \begin{tabular}{ccc}
0.9825 \\
1 \\
1 
\end{tabular}  &  \begin{tabular}{ccc}
    0.969 \\
    1 \\
    1
\end{tabular}  & \begin{tabular}{ccc}
0.955\\
1 \\
1
\end{tabular}\\ 
\hline
\end{tabular}}
\label{Table-power}
\end{table}

From Table \ref{Table-power}, it can be seen that the packing test and the test $\phi^F_n$ work very well against the two heavy-tailed alternatives. Interestingly, the packing test also works well in the chi-squared setting, which is of exponential tail. The asymptotic performances of the Rayleigh test and Bingham test also match with that of Theorem \ref{Rayleigh-asymp} and \ref{Bingham-asymp}. Moreover, unlike the packing test, which exhibits poor empirical type-I error, the proposed test $\phi^F_n$ provides excellent type-I error control while remaining powerful against projected distributions with heavy-tailed marginals.

\section{Discussion} \label{sec-discuss}

In this paper, we have investigated the asymptotic power of three known tests for uniformity on high-dimensional hyperspheres under heavy-tailed alternatives. In particular, we provide the first result on the consistency of the packing test proposed by \cite{Jiang13} against heavy-tailed alternatives. We further show that the two popular tests, the Rayleigh test and the Bingham test, perform poorly against the same class of alternatives. By demonstrating that these three tests are asymptotically independent for large $p$ and $n$, we propose a new test based on Fisher's combination technique that enjoys all the optimality properties of each individual test, and thus is robust against heavy-tailed alternatives while not suffering from the size distortion issue of the packing test.

There are some open questions that could be interesting for future research. Firstly, it is unclear what the asymptotic distribution is when the symmetry assumption in Theorems \ref{Rayleigh-asymp} and \ref{Bingham-asymp} is dropped. In particular, some conditional moment computations become complicated without the symmetry assumption. Secondly, it would be of significant interest to investigate the exact asymptotic distribution of the packing test under heavy-tailed alternatives, as it would allow for asymptotic power calculation. This is equivalent to finding the limiting distribution of the largest off-diagonal entries of the sample correlation matrices. Finally, it would be interesting to develop new procedures that are model-free and of non-parametric nature in the high-dimensional settings. We leave these questions for future research.

\section{Proofs of the main results} \label{sec-proofs}

\subsection{Proof of Theorem \ref{Rayleigh-asymp}}

Our proof is based on martingale central limit theorem (see Corollary 3.1 in \cite{Hall}) and the moment formula in Lemma \ref{heavy moments}. 
Define 
\begin{align*}
    Y_{n,i} :&= \sum_{j=1}^{i-1} \lb \frac{\sqrt{2p}}{n} \bm{X}^{\top}_i \bm{X}_j \rb, \\
        Q_{n} :&= \sum_{i=2}^{n} \mb{E} \lb \lb Y_{n,i} \rb^2 \Big| \bm{X}_1, \bm{X}_2, \cdots, \bm{X}_{i-1} \rb.
\end{align*}
Let us check that $\la Y_{n,i}; 2 \leq i \leq n \ra$ is a martingale difference sequence with respect to the filtration $\sigma_i = \mathcal{F} \lb \bm{X}_1, \bm{X}_2, \cdots, \bm{X}_i \rb$, $2\leq i \leq n$. This is because 
\begin{align*}
    \mb{E} \lb Y_{n,i} \Big|  \bm{X}_1, \bm{X}_2, \cdots, \bm{X}_{i-1}  \rb &= \frac{\sqrt{2p}}{n}\sum_{j=1}^{i-1} \mb{E} \lb \bm{X}^{\top}_i \bm{X}_j \Big| \bm{X}_j \rb = 0
\end{align*}
due to Lemma \ref{heavy-Ustat}. To deduce the asymptotic normality, we need to check the Lindeberg-Feller condition that
\begin{align} \label{lindeberg-cond}
  \sum_{i=2}^{n} \mb{E} \Big[ \lb Y_{n,i} \rb^2 \cdot \mathbf{1}_{\la |Y_{n,i}|>\ve \ra} \Big] \to 0
\end{align}
for every fixed $\ve >0$, and
\begin{align} \label{cond-var}
     Q_{n} \xrightarrow{\mb{P}} 1.
\end{align}
To show \eqref{lindeberg-cond}, it suffices to show that 
\begin{align*}
    \sum_{i=2}^{n} \mb{E} Y_{n,i}^4 \to 0.
\end{align*}
To bound the sum above, observe that 
\begin{align*}
    \mb{E} \Big[  (\bm{X}^{\top}_1 \bm{X}_2 )^2 \cdot  (\bm{X}^{\top}_1 \bm{X}_3 ) \cdot (\bm{X}^{\top}_1 \bm{X}_4 ) \Big] & = \mb{E} \mb{E} \Big[  (\bm{X}^{\top}_1 \bm{X}_2 )^2 \cdot  (\bm{X}^{\top}_1 \bm{X}_3) \cdot (\bm{X}^{\top}_1 \bm{X}_4) 
 \Big| \bm{X}_1 \Big] \\
 &= \mb{E} \Big[ \mb{E} \lb (\bm{X}^{\top}_1 \bm{X}_2 )^{2} \Big| \bm{X}_1 \rb \cdot \mb{E} \lb \bm{X}^{\top}_1 \bm{X}_3  \Big| \bm{X}_1 \rb\\
 & \cdot  \mb{E} \lb \bm{X}^{\top}_1 \bm{X}_4  \Big| \bm{X}_1 \rb \Big] \\
 &= 0,
 \end{align*}
 due to Lemma \ref{heavy-Ustat}. Similarly,
 \begin{align*}
      \mb{E} \Big[  (\bm{X}^{\top}_1 \bm{X}_2)^3 \cdot (\bm{X}^{\top}_1 \bm{X}_3 ) \Big] & = \mb{E} \mb{E} \Big[  (\bm{X}^{\top}_1 \bm{X}_2)^3 \cdot (\bm{X}^{\top}_1 \bm{X}_3 ) \Big| \bm{X}_1 \Big] \\
 & = \mb{E} \Big[ \mb{E} \lb (\bm{X}^{\top}_1 \bm{X}_2 )^3 \Big| \bm{X}_1 \rb \cdot \mb{E} \lb \bm{X}^{\top}_1 \bm{X}_2 \Big| \bm{X}_1 \rb \Big]\\
 & = 0
\end{align*}
again due to Lemma \ref{heavy-Ustat}. Also, 
\begin{align*}
    \mb{E} \Big[  (\bm{X}^{\top}_1 \bm{X}_2 ) \cdot  (\bm{X}^{\top}_1 \bm{X}_3 ) \cdot (\bm{X}^{\top}_1 \bm{X}_4 ) \cdot (\bm{X}^{\top}_1 \bm{X}_5) \Big] 
 = 0,
\end{align*}
which can be seen by conditioning on $\bm{X}_1$ and then use  Lemma \ref{heavy-Ustat}. Thus, we have the bound
\begin{align*}
    \mb{E} \lb Y_{n,i} \rb^4 & \leq \frac{4p^2}{n^4} \cdot \Big[ \sum_{j=1}^{i-1} \mb{E} ( \bm{X}^{\top}_i \bm{X}_j )^4 + C \cdot \sum_{1 \leq r \neq t \leq i-1} \mb{E} \Big[ (\bm{X}^{\top}_i \bm{X}_r )^2 \cdot (\bm{X}^{\top}_i \bm{X}_t )^2 \Big] \Big] \\
    & \leq  \frac{4p^2}{n^4} \cdot \Big[ (i-1) \cdot \mb{E}  (\bm{X}^{\top}_1 \bm{X}_2 )^4 + Ci^2 \cdot \mb{E} \Big[ (\bm{X}^{\top}_1 \bm{X}_2)^2 \cdot (\bm{X}^{\top}_1 \bm{X}_3)^2 \Big] \Big]\\
    &=   \frac{4p^2}{n^4} \cdot \Big[ (i-1) \cdot \frac{\Gamma(2-\alpha/2)}{\Gamma(1-\alpha/2)} \cdot \frac{1}{p}(1+o(1)) + Ci^2 \cdot \frac{1}{p^2} \Big].
\end{align*}
Here the first term in the last line follows from the first statement in Lemma \ref{8th-moment} while the second term follows from conditioning on $\bm{X}_1$ and the second statement in Lemma \ref{heavy-Ustat}. Summing over $2 \leq i \leq n$, we get 
\begin{align*}
    \sum_{i=2}^{n} \mb{E} \lb Y_{n,i} \rb^4 &\leq \frac{4p^2}{n^4} \cdot \frac{\Gamma(2-\alpha/2)}{\Gamma(1-\alpha/2)} \cdot \frac{1}{p}(1+o(1)) \cdot O(n^2)  + \frac{4p^2}{n^4} \cdot \frac{1}{p^2} \cdot O(n^3) \\
    &= O \lb \frac{p}{n^2} \rb + O \lb \frac{1}{n} \rb \to 0
\end{align*}
whenever $p/n^2 \to 0$.  Thus, \eqref{lindeberg-cond} is proven.

Let us now prove \eqref{cond-var}. We will check that 
$$\mb{E} Q_n = 1 + o(1)$$
and 
$$\mbox{Var} \lb Q_n \rb \to 0.$$
The statement regarding the expectation of $Q_n$ is clear since
$$\mb{E} Q_n = \frac{2p}{n^2}  \sum_{1 \leq i<j \leq n} \mb{E} \lb \bm{X}^{\top}_i \bm{X}_j \rb^2 = \frac{2p}{n^2} \cdot \frac{n(n-1)}{2p} = 1 +o(1). $$
Therefore, it remains to bound the variance of $Q_n$. Write
\begin{align*}
            Q_{n} &= \frac{2p}{n^2} \sum_{i=2}^{n} \mb{E} \lb \lb Y_{n,i} \rb^2 \Big| \bm{X}_1, \bm{X}_2, \cdots, \bm{X}_{i-1} \rb \\
   & = \frac{2p}{n^2}\sum_{i=2}^{n} \lb  \sum_{k=1}^{i-1}  \mb{E} \lb (\bm{X}^{\top}_i \bm{X}_k )^{2} \Big| \bm{X}_k \rb + \sum_{1 \leq u \neq v \leq i-1} \mb{E} \lb (\bm{X}^{\top}_i \bm{X}_u ) \cdot (\bm{X}^{\top}_i \bm{X}_v ) \Big| \bm{X}_u, \bm{X}_v \rb \rb \\
   &=  \frac{2p}{n^2} \sum_{i=2}^{n} \lb  \frac{i-1}{p} + \frac{1}{p} \cdot \sum_{1 \leq u \neq v \leq i-1} \bm{X}^{\top}_u \bm{X}^{\top}_v \rb
\end{align*}
where we use Lemma \ref{heavy-Ustat} to get the second term on the last line. Consequently, for some positive integers $a_{ij}$ with $a_{ij} \leq 2n$, we have 
\begin{align*}
\mbox{Var} Q_n &= \frac{4}{n^4} \mbox{Var} \lb \sum_{i=2}^{n} \sum_{1 \leq u \neq v \leq i} \bm{X}^{\top}_u \bm{X}_v \rb \\
& = \frac{4}{n^4} \mbox{Var} \lb \sum_{1 \leq i<j \leq n} a_{ij} \bm{X}^{\top}_i \bm{X}_j \rb \\
& \leq \frac{4}{n^4} \max_{i,j} a^2_{ij} \cdot \frac{n(n-1)}{2} \cdot \mbox{Var} \lb \bm{X}^{\top}_1 \bm{X}_2 \rb = O(p^{-1})
\end{align*}
as $n \to \infty$. Thus \eqref{cond-var} is proved and the proof is completed. $\hfill$ $\square$

\subsection{Proof of Theorem \ref{Bingham-asymp}}

The proof is similar to that of Theorem 1. We again use martingale central limit theorem (see Corollary 3.1 in \cite{Hall}). Define 
\begin{align*}
    \mathcal{Y}_{n,i} &:= \frac{1}{\sqrt{n}} \sum_{1 \leq i<j\leq n} \Big[ \lb \bm{X}^{\top}_i \bm{X}_j\rb^2 - \frac{1}{p} \Big] , \\
    \mathcal{Q}_n &:=  \sum_{i=2}^{n} \mb{E} \lb \lb \mathcal{Y}_{n,i} \rb^2 \Big| \bm{X}_1, \bm{X}_2, \cdots, \bm{X}_{i-1} \rb.
\end{align*}
Lemma \ref{heavy-Ustat} implies that $\la \mathcal{Y}_{n,i}; 2 \leq i \leq n \ra$ is a martingale difference sequence with respect to the sequence of filtration $\sigma_i = \mathcal{F} \lb \bm{X}_1, \bm{X}_2, \cdots, \bm{X}_i \rb$, $2\leq i \leq n$. Moreover, it holds that 
$$\frac{\sqrt{n}}{p} B_n = \sum_{i=2}^{n} \mathcal{Y}_{n,i}$$
where $B_n$ is the Bingham test in \eqref{Bingham}. Thus, to deduce the central limit theorem, it suffices to check the Lyapunov condition
\begin{align} \label{Lyapunov}
        \sum_{i=2}^{n} \mb{E} \mathcal{Y}_{n,i}^4 \to 0.
\end{align}
and 
\begin{align} \label{cond-variance3}
\mathcal{Q}_n \stackrel{\mb{P}}{\to} \frac{(2-\alpha)^2}{8\gamma}
\end{align}
where $\gamma$ is the limit of $p/n$. 

Let us start with the Lyapunov condition \eqref{Lyapunov}. Put $l_n(x)=n^{-1/2} (x^2-p^{-1})$. Note that by Lemma \ref{heavy-Ustat}, we have
\begin{align*}
    \mb{E} \lb l_n(\bm{X}^{\top}_1 \bm{X}_2 ) \Big| \bm{X}_1 \rb &= n^{-1/2} \cdot \mb{E} \lb \mb{E} \lb (\bm{X}^{\top}_1 \bm{X}_2 )^{2} \Big| \bm{X}_1 \rb -\frac{1}{p} \rb =0.
\end{align*}
Thus, we have 
\begin{align*}
    \mb{E} \Big[  l_n(\bm{X}^{\top}_1 \bm{X}_2 )^2 \cdot  l_n(\bm{X}^{\top}_1 \bm{X}_3 ) \cdot l_n(\bm{X}^{\top}_1 \bm{X}_4 ) \Big] & = \mb{E} \mb{E} \Big[  l_n(\bm{X}^{\top}_1 \bm{X}_2 )^2 \cdot  l_n(\bm{X}^{\top}_1 \bm{X}_3) \cdot l_n(\bm{X}^{\top}_1 \bm{X}_4) 
 \Big| \bm{X}_1 \Big] \\
 &= \mb{E} \Big[ \mb{E} \lb l_n(\bm{X}^{\top}_1 \bm{X}_2 )^{2} \Big| \bm{X}_1 \rb \cdot \mb{E} \lb l_n \lb \bm{X}^{\top}_1 \bm{X}_3 \rb \Big| \bm{X}_1 \rb\\
 & \cdot  \mb{E} \lb l_n \lb \bm{X}^{\top}_1 \bm{X}_4 \rb \Big| \bm{X}_1 \rb \Big] \\
 &= 0
 \end{align*}
Similarly,
 \begin{align*}
      \mb{E} \Big[  l_n(\bm{X}^{\top}_1 \bm{X}_2)^3 \cdot l_n(\bm{X}^{\top}_1 \bm{X}_3 ) \Big] & = \mb{E} \mb{E} \Big[  (\bm{X}^{\top}_1 \bm{X}_2)^3 \cdot (\bm{X}^{\top}_1 \bm{X}_3 ) \Big| \bm{X}_1 \Big] \\
 & = \mb{E} \Big[ \mb{E} \lb l_n(\bm{X}^{\top}_1 \bm{X}_2 )^3 \Big| \bm{X}_1 \rb \cdot \mb{E} \lb l_n \lb \bm{X}^{\top}_1 \bm{X}_2 \rb \Big| \bm{X}_1 \rb \Big]\\
 & = 0
\end{align*}
and
\begin{align*}
    \mb{E} \Big[  l_n(\bm{X}^{\top}_1 \bm{X}_2 ) \cdot  l_n(\bm{X}^{\top}_1 \bm{X}_3 ) \cdot l_n(\bm{X}^{\top}_1 \bm{X}_4 ) \cdot l_n(\bm{X}^{\top}_1 \bm{X}_5) \Big] 
 = 0.
\end{align*}
Therefore, for some universal constant $C$, we have the bound
\begin{align*}
    \sum_{i=2}^{n} \mb{E} \mathcal{Y}_{n,i}^4 &\leq \sum_{i=2}^{n} \left[ \sum_{j=1}^{i-1} \mb{E} l_n( \bm{X}^{\top}_i \bm{X}_j )^4 + C \cdot \sum_{1 \leq r \neq t \leq i-1} \mb{E} \Big[ l_n(\bm{X}^{\top}_i \bm{X}_r )^2 \cdot l_n(\bm{X}^{\top}_i \bm{X}_t )^2 \Big] \right] \\
    & \leq \sum_{i=2}^{n}  \left[ (i-1) \cdot \mb{E}  l_n(\bm{X}^{\top}_1 \bm{X}_2 )^4 + Ci^2 \cdot \mb{E} \Big[ l_n(\bm{X}^{\top}_1 \bm{X}_2)^2 \cdot l_n(\bm{X}^{\top}_1 \bm{X}_3)^2 \Big] \right]\\
    & \leq \underbrace{\Big( 1 + 2 + \cdots+ (n-1)  \Big)}_{O(n^2)} \cdot \underbrace{\mb{E}  l_n(\bm{X}^{\top}_1 \bm{X}_2 )^4}_{\text{use $l_n(x)^4 \leq 16n^{-2}\lb x^8 + p^{-4} \rb$ and \eqref{4th-moment-inner}}} \\ 
    &+ C \cdot \underbrace{\Big( 1^2 + 2^2 + \cdots + (n-1)^2 \Big)}_{O(n^3)} \cdot \underbrace{\mb{E} \Big[ l_n(\bm{X}^{\top}_1 \bm{X}_2)^2 \cdot l_n(\bm{X}^{\top}_1 \bm{X}_3)^2 \Big]}_{O(n^{-2}p^{-2}) \ \text{by using \eqref{product-moment}}}  \\
    & \leq   \frac{1}{n^2} \cdot O(n^2) \cdot \left[ \mb{E}(\bm{X}^{\top}_1 \bm{X}_2 )^8 + \frac{1}{p^4}  \right] +  O(n^3) \cdot \frac{1}{n^2p^2}\\
    & \leq \frac{1}{n^2} \cdot  O(n^2) \cdot \left[ \mb{E}(\bm{X}^{\top}_1 \bm{X}_2 )^2 + \frac{1}{p^4}  \right] + O\lb \frac{n}{p^2} \rb \\
    &= O \lb \frac{1}{p} \rb + O \lb \frac{n}{p^2} \rb
\end{align*}
which converges to $0$ since $p/n \to \gamma \in (0,\infty)$. Here we used equations \eqref{4th-moment-inner} and \eqref{product-moment} from Lemma \ref{8th-moment} to deduce the fourth inequality.

It remains to check \eqref{cond-variance3}. Note that 
\begin{align*}
    \mb{E} \mathcal{Q}_n &= \sum_{1 \leq i<j \leq n} \mb{E} l_n^2 \lb \bm{X}^{\top}_i \bm{X}_j \rb \\
    & =  \frac{1}{n} \sum_{1 \leq i<j \leq n } \left[ \mb{E} \lb \bm{X}^{\top}_i \bm{X}_j \rb^4 - \frac{2}{p} \mb{E} \lb \bm{X}^{\top}_i \bm{X}_j \rb^2 + \frac{1}{p^2} \right] \\
    & =  \frac{1}{n} \cdot \frac{n(n-1)}{2} \cdot \left[ \mb{E} \lb \bm{X}^{\top}_1 \bm{X}_2 \rb^4 - \frac{1}{p^2} \right] \\
    & =  \frac{1}{n} \cdot \frac{n(n-1)}{2} \cdot \left[ \frac{1}{p} \cdot \frac{(2-\alpha)^2}{4}\cdot(1+o(1)) - \frac{1}{p^2} \right] \\
    & \to \frac{(2-\alpha)^2}{8\gamma}
\end{align*}
where we use Lemma \ref{8th-moment} in the fourth line. Thus, to finish the proof, we just need to check that $\mbox{Var} \mathcal{Q}_n \to 0$. Note that we have the representation
\begin{align*}
    \mathcal{Q}_n &= \sum_{i=2}^{n-1} \sum_{k=1}^{i-1} \mb{E} \left[ l^2_n \lb \bm{X}^{\top}_i \bm{X}_k \rb \Big | \bm{X}_k \right] +  \sum_{i=2}^{n-1} \sum_{1 \leq u \neq v \leq i-1}  \mb{E} \left[ l_n \lb \bm{X}^{\top}_i \bm{X}_u \rb \cdot l_n \lb \bm{X}^{\top}_i \bm{X}_v \rb \Big | \bm{X}_u, \bm{X}_v \right].
\end{align*}
For $\bm{X}, \bm{Y}, \bm{Z}$ i.i.d. from $\mu_{\alpha,p}$, let us define 
\begin{align}
    h_4(\bm{X}) &:=  \mb{E} \left[  \lb \lb \bm{X}^{\top} \bm{Y} \rb^2- p^{-1} \rb^2 \Big | \bm{X} \right] \label{def-h4},\\
   L(\bm{X}, \bm{Y}) &:=  \mb{E} \left[  \lb \lb \bm{X}^{\top} \bm{Z} \rb^2- p^{-1} \rb \cdot \lb \lb \bm{Y}^{\top} \bm{Z} \rb^2- p^{-1} \rb \Big | \bm{X}, \bm{Y}  \right].  \label{def-L}
\end{align}
One can then simplify $\mathcal{Q}_n$ as
\begin{align*}
    \mathcal{Q}_n &= \frac{1}{n} \sum_{i=2}^{n-1} \sum_{k=1}^{i-1} h_4 \lb \bm{X}_k \rb +  \frac{1}{n} \sum_{i=2}^{n-1}  \sum_{1 \leq u \neq v \leq i-1} L(\bm{X}_u, \bm{X}_v).
\end{align*}
Consequently, we have
\begin{align*}
    \mbox{Var} \mathcal{Q}_n \leq 2 \left[ \mbox{Var}  \lb \frac{1}{n} \sum_{i=2}^{n-1} \sum_{k=1}^{i-1}  h_4 \lb \bm{X}_k \rb \rb + \mbox{Var} \lb \frac{1}{n} \sum_{i=2}^{n-1}  \sum_{1 \leq u \neq v \leq i-1} L(\bm{X}_u, \bm{X}_v) \rb \right] \to 0
\end{align*}
by using Lemma \ref{var bound}. The proof is completed. $\hfill$ $\square$

\subsection{Proof of Theorem \ref{Packing-asymp}}

Let us first define the following $\ve$-good concept, which measures how the largest absolute coordinate influences a vector's length. 

\begin{definition}
    Given a vector $\bm{X}=(X_1, X_2,\cdots, X_p)$, we say that $\bm{X}$ is {\it $\ve$-good} if 
    \begin{align} \label{good}
        \frac{\max_{1\leq i \leq p} |X_i|}{\sqrt{\sum_{i=1}^{p} X_i^2}} \geq 1- \ve.
    \end{align}
\end{definition}
Intuitively, a vector is $\ve$-good if its largest coordinate contribute to at least $(1-\ve)^2$ of its squared Euclidean norm. Equivalently, \eqref{good} says that $|\bm{X}^{\top} \bm{e}_i|$ is large for some vector $\bm{e}_i$ in the canonical basis of $\mb{R}^p$. We also define the mapping 
\begin{align} \label{argmax}
    i(\bm{X})= \argmax_{1 \leq \leq n} |X_i|
\end{align}
which is the index corresponding to the largest coordinate of $\bm{X}$. 

Observe that if $\bm{X}=(X_1, X_2,\cdots,X_n)$ and $\bm{Y}=(Y_1, Y_2,\cdots,Y_n)$ are both $\ve$-good and $i(\bm{X})=i(\bm{Y})=m \in [1,n]$, then 
\begin{align*}
   & \left|  \frac{\sum_{k=1}^{n} X_k Y_k}{ \sqrt{\sum_{k=1}^{n} X_k^2} \cdot \sqrt{\sum_{k=1}^{n} Y_k^2}}  \right| \\
 \geq & \left| \frac{X_m Y_m}{\sqrt{\sum_{k=1}^{n} X_k^2} \cdot \sqrt{\sum_{k=1}^{n} Y_k^2}} \right| - \left| \frac{\sum_{k \neq m} X_k Y_k}{ \sqrt{\sum_{k=1}^{n} X_k^2} \cdot \sqrt{\sum_{k=1}^{n} Y_k^2}} \right| \\
 \geq & (1-\ve)^2 - \left| \sqrt{\frac{\sum_{k \neq m} X_k^2}{\sum_{k=1}^{n} X_k^2}} \cdot \sqrt{\frac{\sum_{k \neq m} Y_k^2}{\sum_{k=1}^{n} Y_k^2}} \right| \\
 \geq & (1-\ve)^2 - \left| \sqrt{1 - \frac{X_m^2}{\sum_{k=1}^{n} X_k^2}} \cdot \sqrt{1 - \frac{Y_m^2}{\sum_{k=1}^{n} Y_k^2}} \right| \\
 \geq & (1-\ve)^2 - \lb 1 - (1-\ve)^2 \rb = 1 - 4\ve + 2\ve^2.
\end{align*}
Based on the observation above, the proof essentially boils down to showing that with high probability, there exist two $\ve$-good vectors $\bm{X}_i$ and $\bm{X}_j$ such that $i(\bm{X}_i)=i(\bm{X}_j)$. 

Fix $\ve>0$ sufficiently small. Recall the positive constant $C_{\alpha,\ve}$ from Lemma \ref{good-prob}. Let $A_n$ be the event that there exist at least $\lfloor n C_{\alpha,\ve} \rfloor/4$ $\ve$-good vectors among $\bm{X}_1, \bm{X}_2,\cdots,\bm{X}_n$. Observe that 
$$A_n =  \la \sum_{i=1}^{n} \mathbbm{1}_{\la \bm{X}_i \text{is $\ve$-good} \ra} \geq \frac{\lfloor n C_{\alpha,\ve} \rfloor}{4} \ra.$$
Since $\bm{X}_i$'s are i.i.d., the indicator functions in the display above are i.i.d. Bernoulli random variables with parameter not smaller than $C_{\alpha,\ve}/2$ (by using Lemma \ref{good-prob}), assuming $n$ is large enough. Thus, by Hoeffding's inequality, we have $\mb{P}(A_n) \geq 1 - e^{-Kn C_{\alpha,\ve}^2}$. Here $K$ is a universal constant. 
Write 
\begin{align*}
   & \mb{P} \lb \max_{1\leq i<j \leq p} |\bm{X}^{\top}_i \bm{X}_j| \geq 1 - 4\ve + 2\ve^2  \rb \\
    \geq & \mb{P} \lb \max_{1\leq i< \leq p} |\bm{X}^{\top}_i \bm{X}_j| \geq 1 - 4\ve + 2\ve^2, A_n  \rb  \\
    = & \mb{P} \lb \max_{1\leq i< \leq p} |\bm{X}^{\top}_i \bm{X}_j| \geq 1 - 4\ve + 2\ve^2 \Big| A_n  \rb \cdot \underbrace{\mb{P}\lb A_n \rb}_{\geq 1 - e^{-Kn C_{\alpha,\ve}^2}} \\
    \geq & \mb{P} \Big( \exists u<v: i(\bm{X}_u)=i(\bm{X}_v) \ \text{and both $\bm{X}_u, \bm{X}_v$ are $\ve$-good} \Big| A_n \Big) \cdot \lb 1 - e^{-Kn C_{\alpha,\ve}^2} \rb.
\end{align*}
The last line in the display above follows from the fact that 
$$\la \exists u<v: i(\bm{X}_u)=i(\bm{X}_v) \ \text{and both $\bm{X}_u, \bm{X}_v$ are $\ve$-good} \ra \subset \la \max_{1\leq i< \leq p} |\bm{X}^{\top}_i \bm{X}_j| \geq 1 - 4\ve + 2\ve^2 \ra.$$
Conditional on the event $A_n$, there exists at least $K \geq \lfloor n C_{\alpha,\ve} \rfloor/4$  vectors $\bm{X}_{i_1}, \cdots,\bm{X}_{i_K}$ which are all $\ve$-good and Lemma \ref{conditional-unif} asserts that $\la i(\bm{X}_{i_j}), 1 \leq j \leq K \ra$ are i.i.d. uniformly  distributed over $\la 1,2,\cdots,p \ra$. Moreover, Lemma \ref{i=j}  yields that,  with probability converging to $1$, we can find $1 \leq u<v \leq p$ such that $i(\bm{X}_u)=i(\bm{X}_v)$ since $p \ll n^2 C_{\alpha,\ve}^2$. Consequently,
\begin{align*}
     \mb{P} \Big( \exists u<v: i(\bm{X}_u)=i(\bm{X}_v) \ \text{and both $\bm{X}_u, \bm{X}_v$ are $\ve$-good} \Big| A_n \Big) \to 1
\end{align*}
as $n \to \infty$ where $\ve>0$ is kept fixed. Thus,
$$\lim_{n \to \infty} \mb{P} \lb \max_{1\leq i< \leq p} |\bm{X}^{\top}_i \bm{X}_j| \geq 1 - 4\ve + 2\ve^2  \rb = 1.$$
The proof is completed since $\ve>0$ is arbitrarily small. $\hfill$ $\square$

\section{Proof of Theorem \ref{asymp-inde}}
Recall $P_n$ from \eqref{packing} and its corresponding asymptotic distribution $G$, which is of Gumbel-type. We need to prove that
\begin{align} \label{asymp-ind}
    \mb{P} \lb R_n \leq x, B_n \leq y, P_n \leq z \rb \to \Phi(x) \cdot \Phi(y) \cdot G(z),
\end{align}
for fixed $(x,y,z) \in \mb{R}^3$.
To prove (\ref{asymp-ind}), it suffices to prove that 
\begin{align}
    \mb{P} \lb R_n \leq x, B_n \leq y \rb &\to \Phi(x) \cdot \Phi(y), \label{asymp1} \\
    \mb{P} \lb P_n \geq z \rb &\to 1- G(z), \label{asymp2}
\end{align}
and 
\begin{align}
    \mb{P} \lb R_n \leq x, B_n \leq y, P_n \geq z \rb - \mb{P} \lb R_n \leq x, B_n \leq y \rb \cdot \mb{P} \lb P_n \geq z \rb \to 0. \label{asymp3}
\end{align}
Note that if (\ref{asymp1}), (\ref{asymp2}) and (\ref{asymp3}) have been granted, then
\begin{align*}
   &  \Big|     \mb{P} \lb R_n \leq x, B_n \leq y, P_n \leq z \rb - \Phi(x) \cdot \Phi(y) \cdot G(z)   \Big| \\
   = & \Big|     \mb{P} \lb R_n \leq x, B_n \leq y \rb -      \mb{P} \lb R_n \leq x, B_n \leq y, P_n \geq z \rb - \Phi(x) \cdot \Phi(y) \cdot G(z)   \Big| \\
   \leq &  \Big|     \mb{P} \lb R_n \leq x, B_n \leq y \rb -  \Phi(x) \cdot \Phi(y)  \Big| + \Big| \mb{P} \lb R_n \leq x, B_n \leq y, P_n \geq z \rb - \Phi(x) \cdot \Phi(y) \cdot (1-G(z)) \Big|   \\ 
   \leq &   \Big|     \mb{P} \lb R_n \leq x, B_n \leq y \rb -  \Phi(x) \cdot \Phi(y)  \Big| \\
   + &  \Big| \mb{P} \lb R_n \leq x, B_n \leq y, P_n \geq z \rb -  \mb{P} \lb R_n \leq x, B_n \leq y \rb \cdot \mb{P} \lb P_n \geq z \rb  \Big| \\
   + & \Big|  \mb{P} \lb R_n \leq x, B_n \leq y \rb \cdot \mb{P} \lb P_n \geq z \rb -  \Phi(x) \cdot \Phi(y) \cdot(1-G(z))  \Big|. 
\end{align*}
The first term and the last term in the expression above tend to $0$ due to \eqref{asymp1} and \eqref{asymp2}, respectively. One can see that the middle term also converges to $0$ by using \eqref{asymp3}. Note that under the assumption $p/ (\log n)^2 \to \infty$, (\ref{asymp2}) follows directly from Theorem 6 in \cite{Jiang13}. In what follows, we will verify (\ref{asymp1}) and (\ref{asymp3}).

The proof of \eqref{asymp1} is based on the martingale central limit theorem and is similar to that of Theorems 1 and 2. The asymptotic independence between $R_n$ and $B_n$ is simply based on the fact that their kernels are uncorrelated and both of them are asymptotically normal. The proof of \eqref{asymp3} is similar to the a classical result in extreme-value theory which states that the partial sum and maximum of a sequence of weakly dependent random variables are asymptotically independent.

 
 {\it \underline{Proof of (\ref{asymp1})}.}  We will show the following joint limit
\begin{align*}
    \begin{pmatrix}
R_n \\
B_n
\end{pmatrix} \stackrel{d}{\to} N\left(\begin{pmatrix}
0 \\
0
\end{pmatrix},\begin{pmatrix}
1 & 0 \\
0 & 1
\end{pmatrix}\right).
\end{align*}
To show the convergence above, it suffices to show that for fixed $(a,b) \in \mb{R}^2$, we have 
$$\sum_{1 \leq i<j \leq n} g_n \lb \bm{X}^{\top}_i \bm{X}_j \rb \to N(0, a^2+b^2),$$
where
\begin{align} \label{def-g}
    g_n(x) = a \cdot \frac{\sqrt{2p}}{n} x + b \cdot \frac{p}{n} \lb x^2- \frac{1}{p} \rb. 
\end{align}
We will proceed via a martingale central limit theorem  (see Corollary 3.1 in \cite{Hall}). Define 
\begin{align*}
        Z_{n1}: &= \sum_{1 \leq i<j \leq n} g_n \lb \bm{X}^{\top}_i \bm{X}_j \rb,  \\
    s_{n1}^2: &= \mbox{Var}(Z_{n1}),\\
    Y^{*}_{n,i} :&= \sum_{j=1}^{i-1} g_n(\bm{X}^{\top}_i \bm{X}_j), \\
        Q_{n1} :&= \sum_{i=2}^{n} \mb{E} \lb \lb Y^{*}_{n,i} \rb^2 \Big| \bm{X}_1, \bm{X}_2, \cdots, \bm{X}_{i-1} \rb.
\end{align*}
We need to check the Lindeberg-Feller condition that
\begin{align} \label{lindeberg-cond2}
   s_{n1}^{-2} \sum_{i=2}^{n} \mb{E} \Big[ \lb Y^{*}_{n,i} \rb^2 \cdot \mathbf{1}_{\la |Y^{*}_{n,i}|>\ve s_{n1} \ra} \Big] \to 0
\end{align}
for every fixed $\ve >0$, and
\begin{align} \label{cond-var2}
    s_{n1}^{-2} Q_{n1} \xrightarrow{\mb{P}} 1.
\end{align}

We split the verifications of \eqref{lindeberg-cond2} and \eqref{cond-var2} into two steps below.

{\it Step 1: Verification the Lindeberg condition (\ref{lindeberg-cond2}).} It suffices to check that
\begin{align} \label{Lindeberg-4th moment}
    s_{n1}^{-4} \sum_{i=2}^{n} \mb{E} \lb Y^{*}_{n,i} \rb^4 \to 0.
\end{align}
Let us show that \eqref{Lindeberg-4th moment} is equivalent to
\begin{align} \label{ratio}
    \frac{\mb{E} \Big[ g_n^4( \bm{X}^{\top}_1 \bm{X}_2 ) \Big]}{n \cdot \mb{E}^2 \Big[  g_n^2(\bm{X}^{\top}_1 \bm{X}_2 ) \Big] } \to 0.
\end{align}
To see this, by using the pairwise independence property (see Lemma \ref{degenerate-Ustat}), we get $s_{n1}^2 = n^2(1+o(1)) \cdot \mb{E} g_n^2(\bm{X}^{\top}_1 \bm{X}_2 )$ and thus, 
\begin{align*} 
    s_{n1}^4 = n^4(1+o(1)) \cdot \lb  \mb{E} g_n^2(\bm{X}^{\top}_1 \bm{X}_2 ) \rb^2 .
\end{align*}
To bound the $4$-th moment terms in (\ref{Lindeberg-4th moment}), we first note that all the mixed moments with at least one odd moment vanish. To see this, write 
\begin{align*}
    \mb{E} \Big[  g_n^2(\bm{X}^{\top}_1 \bm{X}_2 ) \cdot  g_n(\bm{X}^{\top}_1 \bm{X}_3 ) \cdot g_n(\bm{X}^{\top}_1 \bm{X}_4 ) \Big] & = \mb{E} \mb{E} \Big[  g_n^2(\bm{X}^{\top}_1 \bm{X}_2 ) \cdot  g_n(\bm{X}^{\top}_1 \bm{X}_3) \cdot g_n(\bm{X}^{\top}_1 \bm{X}_4) 
 \Big| \bm{X}_1 \Big] \\
 &= \mb{E} \Big[ \mb{E} \lb g_n^2(\bm{X}^{\top}_1 \bm{X}_2 ) \Big| \bm{X}_1 \rb \cdot \mb{E} \lb g_n(\bm{X}^{\top}_1 \bm{X}_3 ) \Big| \bm{X}_1 \rb\\
 & \cdot  \mb{E} \lb g_n(\bm{X}^{\top}_1 \bm{X}_4 ) \Big| \bm{X}_1 \rb \Big] \\
 &= 0,
 \end{align*}
 due to the fact that $\mb{E} g_n (\bm{X}_1^{\top} \bm{X}_2)=0 $ and Lemma \ref{degenerate-Ustat}. Similarly,
 \begin{align*}
      \mb{E} \Big[  g_n^3(\bm{X}^{\top}_1 \bm{X}_2) \cdot  g_n(\bm{X}^{\top}_1 \bm{X}_3 ) \Big] & = \mb{E} \mb{E} \Big[  g_n^3(\bm{X}^{\top}_1 \bm{X}_2) \cdot  g_n(\bm{X}^{\top}_1 \bm{X}_3 ) \Big| \bm{X}_1 \Big] \\
 & = \mb{E} \Big[ \mb{E} \lb g_n^3(\bm{X}^{\top}_1 \bm{X}_2 ) \Big| \bm{X}_1 \rb \cdot \mb{E} \lb g_n(\bm{X}^{\top}_1 \bm{X}_2 ) \Big| \bm{X}_1 \rb \Big]\\
 & = 0
\end{align*}
and 
\begin{align*}
    \mb{E} \Big[  g_n(\bm{X}^{\top}_1 \bm{X}_2 ) \cdot  g_n(\bm{X}^{\top}_1 \bm{X}_3 ) \cdot g_n(\bm{X}^{\top}_1 \bm{X}_4 ) \cdot g_n(\bm{X}^{\top}_1 \bm{X}_5) \Big] 
 = 0,
\end{align*}
which can be seen by conditioning on $\bm{X}_1$ and use  Lemma \ref{degenerate-Ustat}. Consequently, only the terms of the form $\mb{E}X^4$ and $\mb{E}X^2Y^2$ in the expansion of $\mb{E} \lb Y^{*}_{n,i} \rb^4$ are non-zero. Therfore, for some universal constant $C$, we get 
\begin{align*}
    \mb{E} \lb Y^{*}_{n,i} \rb^4 &=  \mb{E}\lb \sum_{j=1}^{i-1} g_n(\bm{X}^{\top}_i \bm{X}_j)\rb^4 \\
    & \leq  \sum_{j=1}^{i-1} \mb{E} g_n^4 ( \bm{X}^{\top}_i \bm{X}_j ) + C \cdot \sum_{1 \leq r \neq t \leq i-1} \mb{E} \Big[ g_n^2(\bm{X}^{\top}_i \bm{X}_r ) \cdot g_n^2(\bm{X}^{\top}_i \bm{X}_t ) \Big]  \\
    &= (i-1) \cdot \mb{E} g_n^4 (\bm{X}^{\top}_1 \bm{X}_2 ) + Ci^2 \cdot \mb{E} \Big[ g_n^2(\bm{X}^{\top}_1 \bm{X}_2) \cdot g_n^2(\bm{X}^{\top}_1 \bm{X}_3) \Big]\\
    &= (i-1) \cdot \mb{E} g_n^4 (\bm{X}^{\top}_1 \bm{X}_2) + Ci^2 \cdot \lb \mb{E} \Big[ g_n^2(\bm{X}^{\top}_1 \bm{X}_2 ) \Big] \rb^2
\end{align*}
where the inequality on the second line follows from the fact that all the terms with odd moments cancel out, and the last line follows from the second statement of Lemma \ref{degenerate-Ustat}. Sum up the above display over all $2 \leq i \leq n$, we arrive at
\begin{align*}
    \sum_{i=2}^{n} \mb{E} \lb Y^{*}_{n,i} \rb^4 &\leq C_1 \cdot \Big(  n^2 \cdot \mb{E} g_n^4 (\bm{X}^{\top}_1 \bm{X}_2) + n^3 \cdot \mb{E}^2 \Big[ g_n^2(\bm{X}^{\top}_1 \bm{X}_2) \Big] \Big) \\
    & \leq (C_1+1) \cdot n^3 \cdot \mb{E} g_n^4 (\bm{X}^{\top}_1 \bm{X}_2),
\end{align*}
for some universal constant $C_1$. This in turn yields
$$s_{n1}^{-4} \sum_{i=2}^{n} \mb{E} \lb Y^{*}_{n,i} \rb^4 \leq 2(C_1+1) \cdot \frac{\mb{E} \Big[ g_n^4( \bm{X}^{\top}_1 \bm{X}_2 ) \Big]}{n \cdot \mb{E}^2 \Big[  g_n^2(\bm{X}^{\top}_1 \bm{X}_2 ) \Big] }.$$
Thus, we only need to prove \eqref{ratio}. To prove \eqref{ratio}, by the Cauchy-Schwartz inequality and a direct computation, we have 
\begin{align*}
    \mb{E} \Big[ g_n^4( \bm{X}^{\top}_1 \bm{X}_2 ) \Big] & \leq 8 \lb  a^4 \cdot \frac{4p^2}{n^4}\mb{E} \lb \bm{X}^{\top}_1 \bm{X}_2   \rb^4 + b^4 \cdot \frac{p^4}{n^4} \mb{E} \left[ (\bm{X}^{\top}_1 \bm{X}_2)^2 - \frac{1}{p}  \right]^4  \rb \\
    & \leq 8 \lb a^4 \cdot \frac{4p^2}{n^4} \cdot \frac{3}{p(p+2)} + b^4 \cdot  \frac{8p^4}{n^4} \left[ \mb{E} \lb \bm{X}^{\top}_1 \bm{X}_2 \rb^8 + \frac{1}{p^4} \right] \rb \\
    & \leq  8 \lb a^4 \cdot \frac{4p^2}{n^4} \cdot \frac{3}{p(p+2)} + b^4 \cdot  \frac{8p^4}{n^4} \left[ \frac{\Gamma(p/2)}{16 \Gamma(4+p/2)} + \frac{1}{p^4} \right] \rb,
\end{align*}
where $\Gamma(x)$ is the standard Gamma function. Here we use the fact that the inner product has explicit moments as described in equations (9.6.2) and (9.6.3) from \cite{M-Jupp}. Note that the term involving Gamma function is asymptotically of order $p^{-4}$ by using the formula $\Gamma(x+1)=x \Gamma(x)$ for all $x>0$. This in turn yields
$$\mb{E} \Big[ g_n^4( \bm{X}^{\top}_1 \bm{X}_2 ) \Big] \leq \frac{C}{n^4}$$
for some universal constant $C>0$. 

Similarly, we have 
\begin{align*}
    \mb{E}^2 \Big[  g_n^2(\bm{X}^{\top}_1 \bm{X}_2 ) \Big] &= \left[ a^2 \cdot \frac{2p}{n^2} \cdot  \frac{1}{p} + b^2 \cdot \frac{p^2}{n^2} \cdot \lb \frac{3}{p(p+2)} - \frac{2}{p^2} + \frac{1}{p^2}    \rb \right]^2 \\
    &= \frac{1}{n^4} \cdot \left[ 2a^2 + b^2 \cdot \lb \frac{3p}{p+2} -1 \rb  \right]\\
    &\geq \frac{2a^2+b^2}{n^4}.
\end{align*}
The two estimates above yield \eqref{ratio}. Consequently, we get (\ref{lindeberg-cond2}).

{\it Step 2: Verification of (\ref{cond-var2})}. To show \eqref{cond-var2}, we only need to show that
\begin{align} \label{H-to-0}
   n^4 \cdot \mb{E} \lb H^2_{n1} \lb \bm{X}, \bm{Y} \rb \rb \to 0
\end{align}
where 

\begin{align*}
    H_{n1} \lb \bm{X}, \bm{Y} \rb &=  \mb{E} \Big[  g_n \lb \bm{X}^{\top} \bm{Z} \rb \cdot g_n \lb  \bm{Y}^{\top} \bm{Z} \rb \Big| \bm{X}, \bm{Y} \Big]
\end{align*}
for i.i.d. $\bm{X}, \bm{Y}, \bm{Z}$ with common distribution $\mbox{Unif} \lb \mb{S}^{p-1} \rb$ and $g_n$ is defined as in (\ref{def-g}). 

To see why \eqref{H-to-0} implies \eqref{cond-var2}, write
\begin{align*}
    Q_{n1} &= \sum_{i=2}^{n} \Big[  \sum_{j=1}^{i-1} \mb{E} \lb g_n^2(\bm{X}^{\top}_i \bm{X}_j) \Big| \bm{X}_j \rb + \sum_{1 \leq r \neq t \leq i-1 } \mb{E} \lb g_n(\bm{X}^{\top}_i \bm{X}_r) \cdot g_n(\bm{X}^{\top}_i \bm{X}_t) \Big| \bm{X}_{r}, \bm{X}_t \rb \Big] \\
    &= \sum_{i=2}^{n} (i-1) \cdot \mb{E} g_n^2(\bm{X}^{\top}_1 \bm{X}_2) + \sum_{i=2}^{n} \sum_{1 \leq r \neq t \leq i-1} \mb{E} \lb g_n(\bm{X}^{\top}_i \bm{X}_r) \cdot g_n(\bm{X}^{\top}_i \bm{X}_t) \Big| \bm{X}_{r}, \bm{X}_t \rb.
\end{align*}
Thus, we have 
\begin{align*}
    Q_{n1} &= \sum_{i=2}^{n} (i-1) \cdot \mb{E} g_n^2(\bm{X}^{\top}_1 \bm{X}_2) + \sum_{i=2}^{n} \sum_{1 \leq r \neq t \leq i-1} H_{n1}(\bm{X}_r, \bm{X}_t)\\
    &= s_{n1}^2 \cdot (1+o(1))  +  Q_n^{*},
\end{align*}
where 
\begin{align*}
    Q_n^{*} &:= \sum_{i=2}^{n} \sum_{1 \leq r \neq t \leq i-1} H_{n1}(\bm{X}_r, \bm{X}_t).
\end{align*}
Thus, proving (\ref{cond-var2}) reduces to showing that $\mbox{Var} \lb Q_n^{*}/s_{n1}^2 \rb \to 0$. Let $A = \la (r,t): 1 \leq r \neq t \leq n-1 \ra$. Note that for any two pair $(r,t) \in A$ and $(r',t') \in A$, we have 
$$\mb{E} \Big[  H_{n1}(\bm{X}_r, \bm{X}_t) \cdot H_{n1}(\bm{X}_{r'}, \bm{X}_{t'}) \Big] =0,$$
unless $\la r,t \ra = \la r', t' \ra$. This is due to the fact that  $\mb{E} g_n (\bm{X}_1^{\top} \bm{X}_2)=0 $. Consequently,
\begin{align*}
    \mbox{Var} \lb Q_n^{*} \rb &= \mbox{Var} \Big( \sum_{(r,t) \in A} \Big[ n - \max \la r,t \ra -1 \Big] \cdot H_{n1}(\bm{X}_r, \bm{X}_t) \Big) \\
    & = \sum_{(r,t) \in A} \Big[ n - \max \la r,t \ra -1 \Big]^2 \cdot \mb{E}  H_{n1}^2(\bm{X}_r, \bm{X}_t) \\
    & \leq O(n^4) \cdot \mb{E}    H_{n1}^2(\bm{X}_1, \bm{X}_2),
\end{align*}
where we use the fact that $|A| \leq n^2$ in the last bound. The last term in the display above is exactly \eqref{H-to-0} and hence, we only need to prove \eqref{H-to-0}.

To show \eqref{H-to-0}, we expand the product inside the expectation to get $    H_{n1} \lb \bm{X}, \bm{Y} \rb= A_1+A_2+A_3+A_4$, where
\begin{align*}
    A_1 &= \frac{2a^2p}{n^2} \cdot \mb{E} \lb \bm{X}^{\top} \bm{Z} \cdot \bm{Y}^\top \bm{Z} \Big| \bm{X}, \bm{Y} \rb,\\
    A_2 &= \frac{\sqrt{2} ab p^{3/2}}{n^2} \cdot \mb{E} \left[ \bm{X}^{\top} \bm{Z} \cdot \lb \lb  \bm{Y}^\top \bm{Z} \rb^2 -\frac{1}{p} \rb \Big| \bm{X}, \bm{Y} \right], \\
    A_3 &=  \frac{\sqrt{2} ab p^{3/2}}{n^2} \cdot \mb{E} \left[ \bm{Y}^{\top} \bm{Z} \cdot \lb \lb  \bm{X}^\top \bm{Z} \rb^2 -\frac{1}{p} \rb \Big| \bm{X}, \bm{Y} \right],\\
    A_4 &=  \frac{b^2 p^{2}}{n^2} \cdot \mb{E} \left[ \lb (\bm{X}^{\top} \bm{Z})^2 - \frac{1}{p} \rb \cdot \lb \lb  \bm{Y}^\top \bm{Z} \rb^2 -\frac{1}{p} \rb \Big| \bm{X}, \bm{Y} \right].
\end{align*}
Easily, we have 
$$A_1 = \frac{2a^2}{n^2} \cdot \bm{X}^{\top} \bm{Y},$$
and

\begin{align*}
    A_2 &=  \frac{\sqrt{2} ab p^{3/2}}{n^2} \cdot \mb{E} \left[ \bm{X}^{\top} \bm{Z} \cdot \lb  \bm{Y}^\top \bm{Z} \rb^2  \Big| \bm{X}, \bm{Y} \right] \\
    &= \frac{\sqrt{2} ab p^{3/2}}{n^2} \cdot  \mb{E} \left[ \bm{X}^{\top} \bm{Z} \cdot \lb  \bm{Y}^\top \bm{Z} \rb \cdot \bm{Z}^{\top} \bm{Y}  \Big| \bm{X}, \bm{Y} \right] \\
    &=  \frac{\sqrt{2} ab p^{3/2}}{n^2} \cdot \bm{X}^{\top} \cdot   \mb{E} \left[  \bm{Z} \cdot \  \bm{Y}^\top \cdot \lb \bm{Z} \bm{Z}^{\top} \rb  \Big| \bm{X}, \bm{Y} \right] \cdot \bm{Y} \\
    &=0,
\end{align*}
where we use the fact that all odd mixed moments of the uniform distribution of hypersphere are $0$ in the last equality. Simiarly, we have $A_3=0$. For the term $A_4$, write

\begin{align*}
    A_4 &= \frac{b^2 p^{2}}{n^2} \cdot \left[ \mb{E} \left[  (\bm{X}^{\top} \bm{Z})^2  \cdot \lb   \bm{Y}^\top \bm{Z} \rb^2  \Big| \bm{X}, \bm{Y} \right] - \mb{E} \left[ \lb   \bm{X}^\top \bm{Z} \rb^2 \Big| \bm{X} \right] - \mb{E} \left[ \lb   \bm{Y}^\top \bm{Z} \rb^2 \Big| \bm{Y} \right] + \frac{1}{p^2} \right] \\
    &= \frac{b^2 p^{2}}{n^2} \cdot \left[ \mb{E} \left[  (\bm{X}^{\top} \bm{Z})^2  \cdot \lb   \bm{Y}^\top \bm{Z} \rb^2  \Big| \bm{X}, \bm{Y} \right] -  \frac{1}{p^2}   \right]. \\
\end{align*}
Let $\bm{X}=(x_1,x_2,\cdots,x_p)$, $\bm{Y}=(y_1, y_2, \cdots, y_p)$ and $\bm{Z}=(z_1, z_2, \cdots, z_p)$. We integrate with respect to $\bm{Z}$ and keep $\bm{X}, \bm{Y}$ fixed. Note that only the even mixed moments of $z_i$'s do not vanish. Write

\begin{align*}
    A_4 &= \frac{b^2 p^{2}}{n^2} \cdot \mb{E}_{z} \left[  \lb \sum_{i=1}^{p} x_iz_i \rb^2 \cdot \lb \sum_{i=1}^{p} y_iz_i \rb^2 - \frac{1}{p^2}  \right]  \\
    &=  \frac{b^2 p^{2}}{n^2} \cdot \mb{E}_{z} \left[  \lb \sum_{i=1}^{p} x^2_i z^2_i + \sum_{j \neq k} x_jz_j x_kz_k \rb \cdot \lb \sum_{l=1}^{p} x^2_l z^2_l + \sum_{s \neq t} x_sz_s x_tz_t \rb - \frac{1}{p^2}  \right].
\end{align*}
Expand the product via the formula (9.6.2) in \cite{M-Jupp} and ignore odd mixed moment terms, we get 
\begin{align*}
    A_4 &=   \frac{b^2 p^{2}}{n^2} \cdot \left[ \frac{3}{p(p+2)} \sum_{i=1}^{p} x^2_i y^2_i + \frac{1}{p(p+2)} \sum_{i \neq j} x^2_i x^2_j + \frac{1}{p(p+2)} \sum_{i \neq j} x_iy_ix_jy_j  - \frac{1}{p^2} \right] \\
    &= \frac{b^2 p^{2}}{n^2} \cdot \left[  \frac{1}{p(p+2)} \sum_{i=1}^{p} x^2_i y^2_i + \frac{1}{p(p+2)} + \frac{\lb \bm{X}^{\top} \bm{Y} \rb^2}{p(p+2)} - \frac{1}{p^2} \right] \\
    &= \frac{b^2 p^{2}}{n^2} \cdot \left[  \frac{1}{p(p+2)} \sum_{i=1}^{p} \lb x^2_i y^2_i -\frac{1}{p^2} \rb + \frac{\lb \bm{X}^{\top} \bm{Y} \rb^2- p^{-1}}{p(p+2)} \right] \\
    &= \frac{b^2 p }{n(p+2)} \cdot \left[  \sum_{i=1}^{p} \lb x^2_i y^2_i -\frac{1}{p^2} \rb + \lb \bm{X}^{\top} \bm{Y} \rb^2- p^{-1} \right],
\end{align*}
where the second equality follows from the fact that $1=(x_1^2+x_2^2+\cdots +x_p^2)(y_1^2+y_2^2+\cdots+y_p^2)$ and $\lb \bm{X}^{\top} \bm{Y} \rb^2 = (x_1y_1+x_2y_2+\cdots+x_py_p)^2=\sum_{i,j} x_iy_ix_jy_j$. Therefore, we deduce that
$$H_{n1} \lb \bm{X}, \bm{Y} \rb = A_1 + A_4 = \frac{2a^2}{n^2} \cdot \bm{X}^{\top} \bm{Y} + \frac{b^2 p }{n(p+2)} \cdot \left[  \sum_{i=1}^{p} \lb x^2_i y^2_i -\frac{1}{p^2} \rb + \lb \bm{X}^{\top} \bm{Y} \rb^2- p^{-1} \right]. $$
By noting that $A_1$ and $A_4$ are uncorrelated, we arrive at
\begin{align*}
    \mb{E} H^2_{n1} \lb \bm{X}, \bm{Y} \rb &= \frac{4a^4}{n^4} \cdot \mb{E} \lb \bm{X}^{\top} \bm{Y} \rb^2 + \frac{b^4}{n^4}(1+o(1)) \cdot \mb{E} \left[ \sum_{i=1}^{p} \lb x^2_i y^2_i -\frac{1}{p^2} \rb + \lb \bm{X}^{\top} \bm{Y} \rb^2- \frac{1}{p}   \right]^2.
\end{align*}
It suffices to check that the last expectation term tends to $0$. To see this, we compute 
\begin{align*}
   & \mb{E} \left[ \sum_{i=1}^{p} \lb x^2_i y^2_i -\frac{1}{p^2} \rb + \lb \bm{X}^{\top} \bm{Y} \rb^2- \frac{1}{p}   \right]^2 \\
    = & \mb{E} \left[  \sum_{i=1}^{p} \lb x^2_i y^2_i -\frac{1}{p^2} \rb \right]^2  + 2\mb{E} \left[  \lb \sum_{i=1}^{p} x^2_iy^2_i - \frac{1}{p} \rb \cdot \lb \lb \bm{X}^{\top} \bm{Y} \rb^2- \frac{1}{p} \rb \right] + \frac{3}{p(p+2)} - \frac{1}{p^2} \\
     = &  \frac{9p}{p^2(p+2)^2} +  \frac{p(p-1)}{p^2(p+2)^2} - \frac{1}{p^2} + 2 \left[ p \cdot \lb \frac{9}{p^2(p+2)^2} + \frac{p-1}{p^2(p+2)^2}\rb - \frac{1}{p^2} \right] + \frac{3}{p(p+2)} - \frac{1}{p^2} \\
     = & \frac{O(1)}{p^2}.
\end{align*}
This completes the proof of (\ref{cond-var2}). Consequently, we get (\ref{asymp1}).

{\it \underline{Proof of (\ref{asymp3})}.} Define 
\begin{align}
    \triangle_n &:= \la (i,j): 1\leq i<j \leq n  \ra, \label{def-triang} \\
    A_n &:= \la R_n \leq x, B_n \leq y \ra, \label{def-An}
\end{align}
and for $I=(i,j) \in \triangle_n$, put
\begin{align}
    C_I = \la p \lb \bm{X}^{\top}_i \bm{X}_j \rb^2 - 4 \log n +\log \log n \geq z \ra. \label{def-CI}
\end{align}
For $I=(a,b) \in \triangle_n$ and $J=(c,d) \in \triangle_n$, we say that $I > J$ if either $a > c$ or $a=c$ and $b>d$. 
Thanks to the inclusion-exclusion principle, for every fixed $k \in \mb{N}$, we have the bounds 
\begin{align}
    \mb{P} \lb \bigcup_{I \in \triangle_n} A_n C_I \rb &\geq \sum_{I_1 \in \triangle_n} \mb{P} \lb A_n C_{I_1} \rb - \sum_{I_1<I_2 }  \mb{P} \lb A_n C_{I_1} C_{I_2} \rb + \cdots \nonumber \\
    & - \sum_{I_1 < I_2 <\cdots<I_{2k}} \mb{P} \lb A_n C_{I_1} C_{I_2} \cdots C_{I_{2k}} \rb. \label{lower-asymp}
\end{align}
and 
\begin{align}
    \mb{P} \lb \bigcup_{I \in \triangle_n} A_n C_I \rb &\leq \sum_{I_1 \in \triangle_n} \mb{P} \lb A_n C_{I_1} \rb - \sum_{I_1<I_2 }  \mb{P} \lb A_n C_{I_1} C_{I_{2}} \rb + \cdots \nonumber \\
    & + \sum_{I_1 < I_2 <\cdots<I_{2k+1}} \mb{P} \lb A_n C_{I_1} C_{I_2} \cdots C_{I_{2k+1}} \rb. \label{upper-asymp}
\end{align}
We will show via \eqref{lower-asymp} and \eqref{upper-asymp} that 
\begin{align} \label{lower2}
     \liminf_{n \to \infty}     \mb{P} \lb \bigcup_{I \in \triangle_n} A_n C_I \rb \geq \liminf_{n \to \infty} \Big[ \mb{P} \lb R_n \leq x, B_n \leq y \rb \cdot \mb{P}\lb P_n \geq z \rb \Big]
\end{align}
and 
\begin{align} \label{upper2}
     \limsup_{n \to \infty}     \mb{P} \lb \bigcup_{I \in \triangle_n} A_n C_I \rb \leq \limsup_{n \to \infty} \Big[ \mb{P} \lb R_n \leq x, B_n \leq y \rb \cdot \mb{P}\lb P_n \geq z \rb \Big].
\end{align}
Given \eqref{lower2} and \eqref{upper2}, the proof is completed by using \eqref{asymp1} and \eqref{asymp2}. We will only show \eqref{lower2} since the proof of \eqref{upper2} is similar. Write
\begin{align}
  &  \sum_{I_1 \in \triangle_n} \mb{P} \lb A_n C_I \rb - \sum_{I_1<I_2 }  \mb{P} \lb A_n C_{I_1} C_{I_2} \rb + \cdots 
    - \sum_{I_1 < I_2 <\cdots<I_{2k}} \mb{P} \lb A_n C_{I_1} C_{I_2} \cdots C_{I_{2k}} \rb \nonumber \\
= & \sum_{I_1 \in \triangle_n} \mb{P} \lb A_n C_I \rb - \sum_{I_1<I_2 }  \mb{P} \lb A_n C_{I_1} C_{I_2} \rb + \cdots 
    + \sum_{I_1 < I_2 <\cdots<I_{2k-1}} \mb{P} \lb A_n C_{I_1} C_{I_2} \cdots C_{I_{2k-1}} \rb \nonumber \\
    - & \sum_{I_1 < I_2 <\cdots<I_{2k}} \mb{P} \lb A_n C_{I_1} C_{I_2} \cdots C_{I_{2k}} \rb \nonumber \\
    = & \mb{P} \lb A_n \rb \left[  \sum_{I_1 \in \triangle_n} \mb{P} \lb C_{I_1} \rb -   \sum_{I_1<I_2 }  \mb{P} \lb C_{I_1} C_{I_2} \rb  + \cdots + \sum_{I_1 < I_2 <\cdots<I_{2k-1}} \mb{P} \lb C_{I_1} C_{I_2} \cdots C_{I_{2k-1}} \rb \right] \label{lower-asymp2} \\
    + & \sum_{j=1}^{2k-1} D(n,j) + N(n,2k)\cdot \mb{P}(A_n), \nonumber
\end{align}
where 
\begin{align}
    D(n,j) &:= \sum_{I_1<I_2<\cdots<I_j} \mb{P} \lb A_n C_{I_1} C_{I_2} \cdots C_{I_{j}} \rb - \mb{P} \lb A_n \rb \cdot \mb{P} \lb C_{I_1} C_{I_2} \cdots C_{I_{j}} \rb, \label{def-Dnj}\\
    N(n,k) &:= \sum_{I_1<I_2<\cdots<I_k} \mb{P} \lb C_{I_1} C_{I_2} \cdots C_{I_{k}} \rb. \label{def-Nnk}
\end{align}
Apply the inclusion-exclusion inequality again in (\ref{lower-asymp2}), we deduce that
\begin{align*}
    & \sum_{I_1 \in \triangle_n} \mb{P} \lb A_n C_I \rb - \sum_{I_1<I_2 }  \mb{P} \lb A_n C_{I_1} C_{I_2} \rb + \cdots 
    - \sum_{I_1 < I_2 <\cdots<I_{2k}} \mb{P} \lb A_n C_{I_1} C_{I_2} \cdots C_{I_{2k}} \rb \\
    \geq & \mb{P} \lb A_n \rb \mb{P} \lb \bigcup_{I \in \triangle_n} C_I \rb + \sum_{j=1}^{2k-1} D(n,j) + N(n,2k)\cdot \mb{P}(A_n).
\end{align*}
By sending $n \to \infty$, we deduce that for every fixed $k$, the right-hand side of (\ref{lower-asymp}) has infimum limit greater than
\begin{align*}
   & \liminf_{n \to \infty} \left[ \mb{P} \lb A_n \rb \mb{P} \lb \bigcup_{I \in \triangle_n} C_I \rb \right] + \sum_{j=1}^{2k-1}\liminf_{n \to \infty} D(n,j) + \liminf_{n \to \infty} \left[  N(n,2k) \mb{P}(A_n) \right] \\
    = &   \liminf_{n \to \infty} \Big[ \mb{P} \lb R_n \leq x, B_n \leq y \rb \cdot \mb{P}\lb P_n \geq z \rb \Big] + \liminf_{n \to \infty} \left[  N(n,2k) \mb{P}(A_n) \right].
\end{align*}
where we used Lemma \ref{dterm} in the last line to get rid of the summand involving $D(n,j)$. Consequently, by first sending $k \to \infty$, then using Lemma \ref{Nterm}, we get \eqref{lower2}. Similarly, we get \eqref{upper2} and this concludes the proof.

\section{Technical results} \label{sec-technical}

\begin{lemma} \label{i=j}
    Let $X_1, X_2, \cdots, X_n$ be independent random variables that are uniformly distributed  over the set $\la 1,2,\cdots,p \ra$. Suppose $p=p_n$ such that $p /n^2 \to 0$. Then, 
    $$ \lim_{n \to \infty} \mb{P} \lb \exists i \neq j: X_i = X_j \rb =  1.$$
\end{lemma}

 {\noindent \textbf{Proof of Lemma \ref{i=j}}.}
   The statement is obvious if $p < n$ due to the pigeonhole principle. When $p=n$, the probability is exactly $1- n! \cdot n^{-n}$, which converges to $1$ as $n \to \infty$. To avoid triviality, we assume $p \geq n+1$. Let $A_n := \la \exists i \neq j: X_i = X_j \ra$, we can write
    \begin{align*}
        \mb{P} \lb A_n \rb &= 1 - \sum_{ (i_1, i_2,\cdots, i_n) \ \text{distinct}} \mb{P} \lb X_1 =i_1, X_2= i_2, \cdots, X_n = i_n \rb \\
        &= 1 - \underbrace{\Big| \la (i_1, i_2,\cdots, i_n) \ \text{such that they are pairwise distinct}  \ra \Big|}_{\frac{p!}{(p-n)!} \text{choices}}\times \frac{1}{p^n} \\
        &= 1 - \frac{1}{p^n} \cdot \frac{p!}{(p-n)!}.
    \end{align*}
   Put $\gamma_n = p/n$.  By using the Stirling approximation formula $n! \sim \sqrt{2 \pi n} \lb n/e \rb^n$, we arrive at
    \begin{align*}
        \mb{P} \lb A_n \rb &= 1 - \frac{1}{p^n} \cdot \sqrt{2 \pi p} \lb \frac{p}{e} \rb^p \cdot \frac{1}{\sqrt{2 \pi(p-n)}} \lb \frac{e}{p-n} \rb^{p-n}\\
        &= 1 - \sqrt{\frac{p}{p-n}} \cdot \exp \la -n \log p - n + p\log p  - (p-n) \log (p-n)  \ra \\
        &= 1 - \sqrt{\frac{\gamma_n}{\gamma_n-1}} \cdot \exp \la  -n \Big(  \log p + 1 - \gamma_n \log p + (\gamma_n -1) \lb \log(1- 1/\gamma_n)+ \log p \rb \Big) \ra \\
        &= 1 - \sqrt{\frac{\gamma_n}{\gamma_n-1}} \cdot \exp \la -n \lb 1 + (\gamma_n -1) \log(1- 1/\gamma_n) \rb \ra.
    \end{align*}
  Consider the function $f(x)=1+(x-1) \log (1-1/x)$ on $(1,\infty)$. It is easy to check that $f'(x)=\log(1-1/x)+1/x \leq 0$ for $x \in (1,\infty)$ so $f$ is decreasing. Moreover, for large $x$, Taylor expansion gives
  $$f(x)= 1 + (x-1)\log (1-1/x)= 1 + (x-1) \lb \frac{-1}{x} + O(x^{-2}) \rb = \frac{1}{x}(1+o(1)).$$
  Therefore, for some universal constant $C>1$, $f(x) \geq 1/2x$ whenever $x \geq C$. Additionally, we also have
  $$f(x) \geq f(C)$$
  for $x \in (1, C]$, due to the decreasing property. Note that $f(C) \geq 1/(2C)>0$. Consequently, if $ \gamma_n \geq C$, then
  \begin{align*}
      \mb{P} \lb A_n \rb &=  1 - \sqrt{\frac{\gamma_n}{\gamma_n-1}} \cdot \exp \la -n \lb 1 + (\gamma_n -1) \log(1- 1/\gamma_n) \rb \ra \\
      &= 1 -   \sqrt{\frac{\gamma_n}{\gamma_n-1}} \cdot \exp \la -n f(\gamma_n) \ra \\
      & \geq 1 - \sqrt{\frac{C}{C-1}} \cdot \exp \lb -C n/ \gamma_n \rb.
  \end{align*}
  On the other hand, for $\gamma_n \in (1+1/n, C)$, we have 
    \begin{align*}
      \mb{P} \lb A_n \rb &=  1 - \sqrt{\frac{\gamma_n}{\gamma_n-1}} \cdot \exp \la -n \lb 1 + (\gamma_n -1) \log(1- 1/\gamma_n) \rb \ra \\
      &= 1 -   \sqrt{\frac{p}{p-n}} \cdot \exp \la -n f(C) \ra \\
      & \geq 1 -  \exp \lb \log n - n \cdot f(C) \rb.
  \end{align*}
  Thus, we have
  \begin{align*}
            \mb{P} \lb A_n \rb & \geq \min \la 1 - \sqrt{\frac{C}{C-1}} \cdot \exp \lb -C n/ \gamma_n \rb,  1 -  \exp \lb \log n - n \cdot f(C) \rb \ra \to 1
  \end{align*}
as $n \to \infty$ since $p/ n^2 \to 0$ and $f(C)>0$. The proof is completed. $\hfill$ $\square$

\begin{lemma} \label{conditional-unif}
    Suppose $\bm{X}_1, \bm{X}_2, \cdots, \bm{X}_n$ are i.i.d. $p$-dimensional random vectors with entries drawn independently from a common distribution $\mu$. Fix $\ve<1 - 2^{-1/2}$. Recall $i(.)$ and the concept $\ve$-good in \eqref{argmax} and \eqref{good}, respectively. Then, conditional on the event that all $\bm{X}_i$'s are $\ve$-good, we have 
    $$i(\bm{X}_i)  \stackrel{\text{i.i.d.}}{\sim} \mbox{Unif} \lb \la 1,2,\cdots,p  \ra \rb.$$
\end{lemma}

 {\noindent \textbf{Proof of Lemma \ref{conditional-unif}}.}
    We need to show that for all integers $a_i$ between $1$ and $p$
    \begin{align*}
        \mb{P} \lb i(\bm{X}_1) =a_1, \cdots, i(\bm{X}_n)=a_n \Big| \forall 1 \leq i \leq n: \bm{X}_i \ \text{is $\ve$-good }\rb = p^{-n}.
    \end{align*}
  It suffices to show that the probability on the LHS above is independent of $a_i$'s.  Write 
    \begin{align*}
              &  \mb{P} \lb i(\bm{X}_1) =a_1, \cdots, i(\bm{X}_n)=a_n \Big| \forall 1 \leq i \leq n: \bm{X}_i \ \text{is $\ve$-good }\rb \\
            = & \frac{\prod_{i=1}^{n} \mb{P} \Big( i(\bm{X}_i)=a_i, \text{ $\bm{X}_i$ is $\ve$-good} \Big)}{ \mb{P} \lb \forall 1 \leq i \leq n: \bm{X}_i \ \text{is $\ve$-good} \rb}   \\
            = & \frac{\prod_{i=1}^{n} \mb{P} \lb  \frac{|X_{i,a_i}|}{\sqrt{\sum_{k=1}^{p} X_{ik}^2}} \geq 1-\ve \rb}{\mb{P} \lb \forall 1 \leq i \leq n: \bm{X}_i \ \text{is $\ve$-good} \rb} \\
            = & \frac{\mb{P} \lb  \frac{|X_{11}|}{\sqrt{\sum_{k=1}^{p} X_{1k}^2}} \geq 1-\ve \rb^n}{\mb{P} \lb \forall 1 \leq i \leq n: \bm{X}_i \ \text{is $\ve$-good} \rb}.
    \end{align*}
    where we have use the distributional invariant of $\bm{X}_i$'s under permutations of coordinates in the last line and the fact that
    $$\la \frac{|X_{i,a_i}|}{\sqrt{\sum_{k=1}^{p} X_{ik}^2}} \geq 1-\ve \ra \subset \la i(\bm{X}_i) =a_i \ra$$
   for all $\ve<1-2^{-1/2}$, in the third line. Since the last probability is independent of $a_i$'s, the proof is completed. $\hfill$ $\square$

\begin{lemma} \label{good-prob}
    Let $\bm{X}=(X_1, X_2, \cdots, X_p)$. Recall the concept of $\ve$-good from \eqref{good}. Suppose $X_i$'s are drawn independently from an $\alpha$-spherical distribution $\mu_{\alpha,p}$ with $\alpha \in (0,2)$. Then,
    $$\lim_{p \to \infty} \mb{P} \Big( \bm{X} \ \text{is $\ve$-good}  \Big) = C_{\alpha,\ve} >0$$
    for all $\ve \in (0,1/2)$. Here $C_{\alpha,\ve} \in (0,1)$ depends only on $\alpha$ and $\ve$. 
\end{lemma}

 {\noindent \textbf{Proof of Lemma \ref{good-prob}}.}
    It is a known result (see \cite{Darling}, Theorem 5.1) that 
    $$\frac{\sum_{i=1}^{p} X_i^2}{\max_{1 \leq i \leq p} X_i^2} \stackrel{d}{\to} \frac{U}{V}$$
as $p \to \infty$. Note that this is true because $X_i^2$'s are i.i.d. and regularly varying with index $\alpha/2 \in (0,1)$. Here the random variables $U,V$ are non-negative, dependent and their ratio has a characteristic function given by
\begin{align} \label{characteristic}
    \mb{E} \exp \lb it \frac{U}{V} \rb = \frac{e^{it}}{1 + \frac{\alpha}{2} \int_{0}^{1} \frac{1- e^{itu}}{u^{\alpha/2 +1}} du}
\end{align}
for $t \in \mb{R}$. 

We split the proof into two steps below. In the first step, we derive an analytic continuation of \eqref{characteristic} to get the the moment generating function (mgf) of $U/V$ and then estimate its tail probability via this mgf.

{\indent \it Step 1: Analytic continuation.} We claim that formula \eqref{characteristic} above extends to the whole upper-half plane $\mb{C}^{+}= \la z \in \mb{C}: \mbox{Im}(z)>0 \ra$. This is done via a version of the maximum modulus principle for the unbounded domain $\mb{C}^{+}$; see Theorem 15.1 in \cite{complex}.

Define
\begin{align*}
    F_1(z) &=     \mb{E} \exp \lb iz \frac{U}{V} \rb, \\
    F_2(z) &= \frac{e^{iz}}{1 + \frac{\alpha}{2} \int_{0}^{1} \frac{1- e^{izu}}{u^{\alpha/2 +1}} du}.
\end{align*}
To deduce that $F_1(z)=F_2(z)$ for all $z \in \mb{C}^{+}$ via the maximum modulus principle (Theorem 15.1 in \cite{complex}), we need to check that (i) both $F_1$ and $F_2$ are analytic in $\mb{C}^{+}$ and (ii) $|F_1(z)-F_2(z)|$ are uniformly bounded in $\mb{C}^{+}$.

Let us start with (i). It is easy to see that $F_1$ is analytic since with $z=a+bi$, we have
\begin{align*}
    \frac{d}{dz} F_1(z) &=  i \cdot \mb{E} \lb  \frac{U}{V} \exp \lb iz \frac{U}{V} \rb \rb \\
    &= i \cdot \mb{E} \lb  \frac{U}{V} \exp \lb ia\frac{U}{V} - b \frac{U}{V} \rb \rb.
\end{align*}
The last term converges absolutely since $b=\mbox{Im}(z)>0$. 

To check that $F_2$ is analytic, we only need to check that its denominator is analytic and is non-zero over $\mb{C}^{+}$. Let us start with the former statement, pick a simple closed $C^1$ curve $\tau \subset \mb{C}^{+}$ and write 
\begin{align*}
     & \int_{\tau} \left[ 1 + \frac{\alpha}{2} \int_{0}^{1} \frac{1- e^{izu}}{u^{\alpha/2 +1}} du \right] dz \\
     = & \int_{\tau} 1 dz + \frac{\alpha}{2}  \cdot \int_{\tau} \left[ \int_{0}^{1} \frac{1- e^{izu}}{u^{\alpha/2 +1}} du \right] dz \\
     = &0 +  \int_{0}^{1}  \left[ \int_{\tau} \frac{1- e^{izu}}{u^{\alpha/2 +1}} dz \right] du \\
     = & 0,
\end{align*}
where we used the fact that for every $u \in [0,1]$, $z \mapsto e^{izu}$ is an analytic function on $C^{+}$, in the last line. Thus, by Morera's theorem, the denominator of $F_2$ is an analytic function on $\mb{C}^{+}$. 

To see why the denominator of $F_2$ is nonzero, write 
$$\left|  1 + \frac{\alpha}{2} \int_{0}^{1} \frac{1- e^{izu}}{u^{\alpha/2 +1}} du  \right| \geq \mbox{Im} \lb 1 + \frac{\alpha}{2} \int_{0}^{1} \frac{1- e^{izu}}{u^{\alpha/2 +1}} du \rb = 1 + \frac{\alpha}{2} \int_{0}^{1} \frac{1- e^{-bu}\cdot \cos(au)}{u^{\alpha/2 +1}} du \geq 1$$
where $z=a+bi$ in the display above. Thus, both $F_1$ and $F_2$ are well-defined, analytic function on $\mb{C}^{+}$. Let us now check the boundedness assumption (ii). Write 
\begin{align*}
    |F_1(z) - F_2(z)| &= \left| \mb{E} \exp \lb iz \frac{U}{V} \rb - \frac{e^{iz}}{1 + \frac{\alpha}{2} \int_{0}^{1} \frac{1- e^{izu}}{u^{\alpha/2 +1}} du}  \right| \\
    & \leq 1 + \frac{\exp \lb -\mbox{Im}(z) \rb}{\mbox{Im} \lb 1 + \frac{\alpha}{2} \int_{0}^{1} \frac{1- e^{izu}}{u^{\alpha/2 +1}} du \rb} \\
    & \leq 1 + \frac{\exp \lb -\mbox{Im}(z) \rb}{1 + \frac{\alpha}{2} \int_{0}^{1} \frac{1- e^{-bu}\cdot \cos(au)}{u^{\alpha/2 +1}} du} \\
   & \leq 1 + \exp \lb -\mbox{Im}(z) \rb \leq 2.
\end{align*}
Therefore, (ii) is proven. Hence, $F_1(z)=F_2(z)$ for all $z\in \mb{C}^{+}$ and by restricting to the ray $ \la \lambda i: \lambda>0 \ra$, we obtain 
\begin{align} \label{mgf-bound}
    \mb{E} \exp \lb -\lambda \frac{U}{V} \rb = \frac{e^{-\lambda}}{1 + \frac{\alpha}{2} \int_{0}^{1} \frac{1 - e^{-\lambda u}}{u^{\alpha/2+1}} du}
\end{align}
for all $\lambda>0$.

{\indent \it Step 2: Tail bound.} By definition, we can see that 
    $$\lim_{p \to \infty} \mb{P} \Big( \bm{X} \ \text{is $\ve$-good}  \Big) = \mb{P} \lb \frac{V}{U} \geq (1-\ve)^2 \rb.$$
It suffices to prove that the last probability is strictly positive. Assume the contrary, which means $V/U \leq (1-\ve)^2$ almost surely. Then, we have
\begin{align*}
    \mb{E} \lb -\lambda \frac{U}{V} \rb \leq \exp \lb \frac{-\lambda}{(1-\ve)^2} \rb 
\end{align*}
for all $\lambda > 0$. In the view of \eqref{mgf-bound}, this yields the bound
\begin{align*}
   & \frac{e^{-\lambda}}{1 + \frac{\alpha}{2} \int_{0}^{1} \frac{1 - e^{-\lambda u}}{u^{\alpha/2+1}} du} \leq \exp \lb \frac{-\lambda}{(1-\ve)^2} \rb \\
   \Leftrightarrow & \exp \left[ \lambda \lb \frac{1}{(1-\ve)^2}-1 \rb \right] \leq 1+ \frac{\alpha}{2}\int_{0}^{1} \frac{1 - e^{-\lambda u}}{u^{\alpha/2+1}} du \\
   \Leftrightarrow & \exp \left[ \frac{\lambda\ve(2-\ve)}{(1-\ve)^2} \right] \leq e^{-\lambda} + \lambda \int_{0}^{1} \frac{e^{-\lambda u}}{u^{\alpha/2}} du.
\end{align*}
By letting $\lambda \to \infty$, we get a contradiction since the LHS grows exponentially fast in $\lambda$ while the RHS is of order $O(\lambda)$. The proof is completed. $\hfill$ $\square$

\begin{lemma}[\cite{heiny2022}] \label{heavy moments}
    Let $\bm{X}=\lb X_1, X_2,\cdots, X_p \rb$ with the symmetric $\alpha$-spherical distribution $\mu_{\alpha,p}$ for some $\alpha \in (0,2)$. Then, for all $1 \leq r \leq p$, we have 
    \begin{align} \label{heavy-mixed-monent}
     \binom{p}{r} \mb{E} \lb X_{i_1}^{2k_1} X_{i_2}^{2k_2} \cdots X_{i_r}^{2k_r} \rb \sim \frac{\left(\frac{\alpha}{2}\right)^{r-1} \prod_{j=1}^r \Gamma\left(k_j-\alpha / 2\right)}{r\Big( \Gamma(1-\alpha / 2) \Big)^r \Gamma(k)}
    \end{align}
    as $p \to \infty$. 
\end{lemma}

 {\noindent \textbf{Proof of Lemma \ref{heavy moments}}.}
    The result follows from equation (3.2) in \cite{heiny2022} and the fact that the distribution of $\bm{X}$ is invariant under permutation.  $\hfill$ $\square$

\begin{lemma} \label{heavy-Ustat}
    Let $\bm{X}_i=\lb X_{i1}, X_{i2},\cdots, X_{ip} \rb$, $1\leq i \leq n$ be i.i.d. with a symmetric $\alpha$-spherical distribution law $\mu_{\alpha,p}$. Here  $\alpha \in (0,2)$. Then, we have 
    \begin{enumerate}
        \item $ \mb{E} \lb \bm{X}_1^{\top} \bm{X}_2 \Big| \bm{X}_1 \rb = 0$ almost surely.
        \item $\mb{E} \left[ \lb \bm{X}_1^{\top} \bm{X}_2 \rb^2 \Big| \bm{X}_1 \right] = 1/p$ almost surely.
        \item $\mb{E} \left[ \lb \bm{X}_1^{\top} \bm{X}_2 \rb \cdot \lb \bm{X}_1^{\top} \bm{X}_3 \rb  \Big| \bm{X}_2, \bm{X}_3 \right] = p^{-1} \cdot \bm{X}^{\top}_2 \bm{X}_3$ almost surely.
    \end{enumerate}
\end{lemma}

 {\noindent \textbf{Proof of Lemma \ref{heavy-Ustat}}.}
    For the first claim, write 
    \begin{align*}
         \mb{E} \lb \bm{X}_1^{\top} \bm{X}_2 \Big| \bm{X}_1 \rb &= \sum_{k=1}^{p} X_{1k} \mb{E} \lb X_{2k} \Big| \bm{X}_1 \rb \\
         &= \sum_{k=1}^{p} X_{1k} \mb{E} \lb X_{2k} \rb = 0.
    \end{align*}
    The second equality in the expression below holds almost surely so the first statement follows. To prove the second statement, write 
    \begin{align*}
         \mb{E} \Big[ \lb \bm{X}_1^{\top} \bm{X}_2 \rb^2 \Big| \bm{X}_1 \Big] &= \sum_{k=1}^{p} X_{1k}^2 \mb{E} \lb X_{2k}^2 \Big| \bm{X}_1 \rb + 2\sum_{1 \leq i <j \leq p} X_{1i} X_{1j} \mb{E} \lb X_{2i} X_{2j} \Big| \bm{X}_1 \rb \\
         &= p^{-1} \sum_{k=1}^{p} X_{1k}^2 + 2\sum_{1 \leq i <j \leq p} X_{1i} X_{1j} \mb{E} \lb X_{2i} X_{2j} \rb =p^{-1}
    \end{align*}
    almost surely, where the last equality follows from the symmetry assumption, which gives $\mb{E} X_{2i} X_{2j} =0$ for all $i \neq j$.
    For the third statement, write 
    \begin{align*}
        \mb{E} \left[ \lb \bm{X}_1^{\top} \bm{X}_2 \rb \cdot \lb \bm{X}_1^{\top} \bm{X}_3 \rb  \Big| \bm{X}_2, \bm{X}_3 \right]  &= \mb{E} \left[ \lb \sum_{k=1}^{p} X_{1k} X_{2k}\rb \cdot \lb \sum_{k=1}^{p} X_{1k} X_{3k}\rb  \Big| \bm{X}_2, \bm{X}_3 \right]  \\
        &= \sum_{k=1}^{p} X_{2k} X_{3k} \cdot \mb{E}X_{1k}^2 + \sum_{1 \leq i \neq j \leq p} X_{2i} X_{3j} \mb{E} X_{1i} X_{1j} \\
        &= p^{-1} \cdot \bm{X}^{\top}_2 \bm{X}_3 +0 =  p^{-1} \cdot \bm{X}^{\top}_2 \bm{X}_3
    \end{align*}
    almost surely. $\hfill$ $\square$

\begin{lemma} \label{8th-moment}
    With the same notation as in Lemma \ref{heavy-Ustat}, as $p \to \infty$, we have
    \begin{align} \label{4th-moment-inner}
            \mb{E} \lb \bm{X}^{\top}_1 \bm{X}_2 \rb^4  = \frac{1}{p} \cdot \lb \frac{\Gamma(2-\alpha/2)}{\Gamma(1-\alpha/2)} \rb^2(1+o(1)) &= \frac{1}{p} \cdot \lb \frac{2-\alpha}{2} \rb^2(1+o(1))
    \end{align}
    and
    \begin{align} \label{product-moment}
          \mb{E} \left[ \lb \lb \bm{X}_1^{\top} \bm{X}_2 \rb^2 - p^{-1} \rb^2 \cdot  \lb \lb \bm{X}_1^{\top} \bm{X}_3 \rb^2 - p^{-1} \rb^2  \right] &= O(p^{-2}).
    \end{align}
\end{lemma}

 {\noindent \textbf{Proof of Lemma \ref{8th-moment}}.}
Let us start with the second statement. Observe that
\begin{align*}
   & \mb{E} \left[ \lb \bm{X}_1^{\top} \bm{X}_2 \rb^2 \lb \bm{X}_1^{\top} \bm{X}_3 \rb^2 \right]  \\
   = & \mb{E} \left[  \mb{E} \left[ \lb \bm{X}_1^{\top} \bm{X}_2 \rb^2 \Big| \bm{X}_1 \right] \cdot \mb{E} \left[ \lb \bm{X}_1^{\top} \bm{X}_3 \rb^2 \Big| \bm{X}_1 \right] \right] \\
   = & \frac{1}{p^2}. 
\end{align*}
Thus, we get
\begin{align*}
    &  \mb{E} \left[ \lb \lb \bm{X}_1^{\top} \bm{X}_2 \rb^2 - p^{-1} \rb^2 \cdot  \lb \lb \bm{X}_1^{\top} \bm{X}_3 \rb^2 - p^{-1} \rb^2  \right]  \\
    = & \mb{E} \left[ \lb \lb \bm{X}_1^{\top} \bm{X}_2 \rb^4 - \frac{2 \lb \bm{X}_1^{\top} \bm{X}_2 \rb^2}{p} +\frac{1}{p^2} \rb \cdot \lb \lb \bm{X}_1^{\top} \bm{X}_3 \rb^4 - \frac{2 \lb \bm{X}_1^{\top} \bm{X}_3 \rb^2}{p} +\frac{1}{p^2} \rb \right] \\
    = & \mb{E} \left[ \lb \bm{X}_1^{\top} \bm{X}_2 \rb^4 \lb \bm{X}_1^{\top} \bm{X}_3 \rb^4 \right] -\frac{2 \mb{E} \left[  \lb \bm{X}_1^{\top} \bm{X}_2 \rb^4 \lb \bm{X}_1^{\top} \bm{X}_3 \rb^2 + \lb \bm{X}_1^{\top} \bm{X}_2 \rb^2 \lb \bm{X}_1^{\top} \bm{X}_3 \rb^4  \right]}{p} +O(p^{-2}) \\
    \leq & \mb{E} \left[ \lb \bm{X}_1^{\top} \bm{X}_2 \rb^2 \lb \bm{X}_1^{\top} \bm{X}_3 \rb^2 \right] + O(p^{-2}) = O(p^{-2}).
\end{align*}
Now we prove the first statement. Expand the fourth moment and note that all the terms involving odd orders cancel out due to symmetry, we arrive at
\begin{align*}
    \mb{E} \lb \bm{X}^{\top}_1 \bm{X}_2 \rb^4 &= \mb{E} \lb \sum_{k=1}^{p} X_{1k}X_{2k} \rb^4 \\
    &= \mb{E} \lb \sum_{k=1}^{p} X_{1k}^4X_{2k}^4  \rb + 2\mb{E} \lb \sum_{1 \leq i<j \leq p} X^2_{1i} X^2_{1j} X^2_{2i} X^2_{2j}\rb \\
    &= p \lb \mb{E} X_{11}^4 \rb^2 + 2 \cdot \frac{p(p-1)}{2} \lb \mb{E} X_{12}^2X_{13}^2 \rb^2.
\end{align*}
Thanks to Lemma \ref{heavy moments}, we can simplify the expression above as
\begin{align*}
        \mb{E} \lb \bm{X}^{\top}_1 \bm{X}_2 \rb^4 &= p \cdot \frac{1}{p^2} \cdot \lb \frac{\Gamma(2-\alpha/2)}{\Gamma(1-\alpha/2)\Gamma(2)} \rb^2 + O(p^2) \cdot \lb \frac{2}{p(p-1)} \cdot \frac{\alpha}{2} \cdot \frac{(\alpha/2)\Gamma(2-\alpha)^2}{2 \Gamma(1-\alpha/2)^2 \Gamma(4)} \rb^2 \\
        &= \frac{1}{p} \cdot \underbrace{\lb \frac{\Gamma(2-\alpha/2)}{\Gamma(1-\alpha/2)} \rb^2}_{(1-\alpha/2)^2} + O \lb \frac{1}{p^2} \rb \\
        &= \frac{1}{p} \cdot \lb \frac{2-\alpha}{2} \rb^2(1+o(1)).
\end{align*}
The proof is completed. $\hfill$ $\square$

\begin{lemma} \label{var bound}
    Recall $h_4$ and $L$ from \eqref{def-h4} and \eqref{def-L}, respectively. Let $\bm{X}_i$'s be the same as in Lemma \ref{heavy-Ustat}. Then, we have 
    $$\lim_{n \to \infty} \mbox{Var}  \lb \frac{1}{n} \sum_{i=2}^{n-1} \sum_{k=1}^{i-1}  h_4 \lb \bm{X}_k \rb \rb = 0 $$\
    and
    $$\lim_{n \to \infty} \mbox{Var} \lb \frac{1}{n} \sum_{i=2}^{n-1}  \sum_{1 \leq u \neq v \leq i-1} L(\bm{X}_u, \bm{X}_v) \rb =0.$$
\end{lemma}

 {\noindent \textbf{Proof of Lemma \ref{var bound}}.}
    Let us start with the first statement. Note that since $h_4(\bm{X}_i)$'s are i.i.d., we have
    \begin{align*}
        & \mbox{Var}  \lb \frac{1}{n} \sum_{i=2}^{n-1} \sum_{k=1}^{i-1}  h_4 \lb \bm{X}_k \rb \rb \\
        = & \frac{1}{n^2} \mbox{Var} \lb \sum_{i=1}^{n-1} (n-i) \cdot h_4 (\bm{X}_i) \rb \\
        = & \frac{1}{n^2}  \cdot O(n^3) \cdot \mbox{Var} \lb h_4 (\bm{X}_1) \rb.
    \end{align*}
    It suffices to show that $n \cdot \mbox{Var} \lb h_4 (\bm{X}_1) \rb \to 0 $. Observe that
    \begin{align*}
      0 \leq  h_4 (\bm{X}_1) &=  \mb{E} \left[ \lb \lb \bm{X}^{\top}_1 \bm{X}_2 \rb^2 - \frac{1}{p} \rb^2 \Big| \bm{X}_1\right] \\
       &=   \mb{E} \left[ \lb \bm{X}^{\top}_1 \bm{X}_2 \rb^4 \Big| \bm{X}_1 \right] - \frac{1}{p^2} \\
       & \leq    \mb{E} \left[ \lb \bm{X}^{\top}_1 \bm{X}_2 \rb^2 \Big| \bm{X}_1 \right] - \frac{1}{p^2} = \frac{1}{p} - \frac{1}{p^2}.
    \end{align*}
    Thus, 
    $$n \cdot \mbox{Var} \lb h_4 (\bm{X}_1) \rb \leq n \mb{E} h^2_4 \lb \bm{X}_1 \rb = O \lb \frac{n}{p^2} \rb$$
    which converges to $0$ since $p/n \to \gamma >0$.

    It remains to show the second claim. To see this, note that   we have the expansion  $\sum_{i=2}^{n-1}  \sum_{1 \leq u \neq v \leq i-1} L(\bm{X}_u, \bm{X}_v) = \sum_{1 \leq u<v \leq n} a_{uv} \cdot L(\bm{X}_u, \bm{X}_v)$
    for some positive integers $a_{uv}$ such that $\max_{u,v} a_{uv} \leq 2n$. Moreover, $L(\bm{X}_u, \bm{X}_v) $ and $L(\bm{X}_s, \bm{X}_t) $ are uncorrelated, unless $\la u,v \ra =  \la s,t \ra$. Therefore, we have 
    \begin{align*}
        & \mbox{Var} \lb \frac{1}{n} \sum_{i=2}^{n-1}  \sum_{1 \leq u \neq v \leq i-1} L(\bm{X}_u, \bm{X}_v) \rb \\
        = & \frac{1}{n^2} \sum_{1 \leq u<v \leq n} a_{uv}^2 \cdot \mbox{Var} \lb  L(\bm{X}_u, \bm{X}_v) \rb \\
         \leq & \frac{1}{n^2} \cdot \lb \max_{u,v} a^2_{uv} \rb \cdot \frac{n(n-1)}{2} \cdot  \mbox{Var} \lb  L(\bm{X}_1, \bm{X}_2) \rb.
    \end{align*}
To finish the proof, it suffices to show that $n^2 \cdot \mbox{Var} \lb  L(\bm{X}_1, \bm{X}_2) \rb \to 0$. By definition, we have
\begin{align*}
L \lb \bm{X}_1, \bm{X}_2) \rb &=  \mb{E} \left[  \lb \lb \bm{X}_1^{\top} \bm{X}_3 \rb^2- p^{-1} \rb \cdot \lb \lb \bm{X}^{\top}_2 \bm{X}_3 \rb^2- p^{-1} \rb \Big | \bm{X}_1, \bm{X}_2  \right] \\
&= \mb{E} \left[ \lb \bm{X}_1^{\top} \bm{X}_3 \rb^2 \cdot \lb \bm{X}_2^{\top} \bm{X}_3 \rb^2  \Big | \bm{X}_1, \bm{X}_2 \right] - \frac{1}{p^2} \\
&= \mb{E} \left[ \lb \sum_{k=1}^{p} X_{1k} X_{3k} \rb^2 \cdot \lb \sum_{k=1}^{p} X_{2k} X_{3k} \rb^2 \Big | \bm{X}_1, \bm{X}_2 \right]  - \frac{1}{p^2} \\
&= \lb \mb{E} X_{31}^4 \rb \cdot \sum_{k=1}^{p} X_{1k}^2 X_{2k}^2 + \lb \mb{E} X_{31}^2 X_{32}^2 \rb \cdot \sum_{1 \leq i \neq j \leq p} X_{1i}^2 X_{2j}^2 - \frac{1}{p^2} \\
&= O \lb \frac{1}{p} \rb \cdot \sum_{k=1}^{p} X_{1k}^2 X_{2k}^2  + O \lb \frac{1}{p^2}\rb
\end{align*}
where we used Lemma \ref{heavy moments} to deduce that $\mb{E} X_{31}^2 X_{32}^2 = O(p^{-2})$ and $\mb{E} X_{31}^4 = O(p^{-1})$ in the last line. Consequently,
\begin{align*}
    n^2 \cdot  \mbox{Var} \lb  L(\bm{X}_1, \bm{X}_2) \rb &= n^2 \cdot \mb{E} L^2(\bm{X}_1, \bm{X}_2)  \\
    &\leq \frac{n^2}{p^2} \cdot \mb{E} \lb  \sum_{k=1}^{p} X_{1k}^2 X_{2k}^2  \rb^2 + O \lb \frac{n^2}{p^4} \rb \\
    &= \frac{n^2}{p^2} \cdot \mb{E} \left[ \sum_{k=1}^{p} X_{1k}^4 X_{2k}^4 + 2 \sum_{1 \leq i<j \leq p} X_{1i}^2 X_{1j}^2 X_{2i}^2 X_{2j}^2 \right] +  O \lb \frac{n^2}{p^4} \rb \\
    & = \frac{n^2}{p^2} \cdot p \cdot \mb{E} \lb X_{11}^4\rb^2 +  \frac{n^2}{p^2} \cdot \frac{p(p-1)}{2} \cdot \mb{E} \left[ X_{11}^2 X_{12}^2 \right]^2 + O \lb \frac{n^2}{p^4} \rb\\
    &= O \lb \frac{n^2}{p^3} \rb + O \lb \frac{n^2}{p^4} \rb.
\end{align*}
where we again use  Lemma \ref{heavy moments} to deduce that $\mb{E} X_{11}^2 X_{12}^2 = O(p^{-2})$ and $\mb{E} X_{11}^4 = O(p^{-1})$ in the last line. Since $p$ and $n$ are proportional, the last display tends to $0$. The proof is completed. $\hfill$ $\square$

\begin{lemma} \label{degenerate-Ustat}
Let $\bm{X}$, $\bm{Y}$ and $\bm{Z}$ be i.i.d realizations of the uniform distribution on $\mb{S}^{p-1}$ and $f: \mb{R} \mapsto \mb{R}$ be a measurable function such that $\mb{E} |f(\bm{X}^{\top} \bm{e}_1)|< \infty$, where $\bm{e}_1= (1,0,\cdots,0)$. Then, we have 
   \begin{itemize}
       \item 
    $\mb{E} \lb f(\bm{X}^{\top} \bm{Y}) \Big|  \bm{Y} \rb = \mb{E} f(\bm{X}^{\top} \bm{Y}) $
    almost surely. 
    \item $\bm{X}^{\top} \bm{Y}$ and $\bm{X}^{\top} \bm{Z}$ are independent.
   \end{itemize}
    
\end{lemma}

 {\noindent \textbf{Proof of Lemma \ref{degenerate-Ustat}}.}
    The first claim is a consequence of the rotational invariant property. Conditioning on $\bm{Y}$, there exists an orthogonal matrix $O$ such that $O^{\top} \bm{Y} = \bm{e}_1= (1,0,\cdots,0)$. Thus, with probability one, we have
\begin{align*}
    \mb{E} \lb f(\bm{X}^{\top} \bm{Y}) \Big|  \bm{Y} \rb &= \mb{E} \lb f \lb \bm{X}^{\top} O^{\top} \bm{Y} \rb \Big|  \bm{Y} \rb \\
    &= \mb{E} \lb f(\bm{X}^{\top} \bm{e}_1) \Big|  \bm{Y} \rb\\
    &= \mb{E} \lb f(\bm{X^{\top}} \bm{e}_1) \rb ,
\end{align*}
where we use the fact that $\bm{X}^{\top} \bm{e}_1$ is independent from $\bm{Y}$ in the last equality. The conditional expectation is well-defined since  Similarly, one can also show that $\mb{E} f(\bm{X}^{\top} \bm{Y}) = \mb{E} \lb f(\bm{X^{\top}} \bm{e}_1) \rb$. This concludes the proof of the first claim.

For the second claim, take any bounded measurable functions $f$ and $g$, by conditioning on $\bm{X}$, we get
\begin{align*}
    \mb{E} \lb f(\bm{X}^{\top} \bm{Y}) g(\bm{X}^{\top} \bm{Z}) \rb &= \mb{E} \Big[ \mb{E} \lb f(\bm{X}^{\top} \bm{Y}) \Big| \bm{X} \rb \cdot \mb{E} \lb g(\bm{X}^{\top} \bm{Z}) \Big| \bm{X} \rb  \Big] \\
    &= \mb{E} f(\bm{X}^{\top} \bm{Y}) \cdot \mb{E} g(\bm{X}^{\top} \bm{Z}),
\end{align*}
where we use the conclusion of the first statement in the last equality. This concludes the proof. $\hfill$ $\square$

\begin{lemma} \label{Nterm}
Recall $\triangle_n$, $A_n$, $C_I$, $N(n,k)$ defined in (\ref{def-triang}), (\ref{def-An}), (\ref{def-CI}), and (\ref{def-Nnk}) respectively. Then, for every fixed $k \in \mb{N}$, we have 
$$\lim_{k \to \infty} \limsup_{n \to \infty} N(n,k) = 0.$$
\end{lemma}

 {\noindent \textbf{Proof of Lemma \ref{Nterm}}.}
    This is essentially the content of Lemma 35 in \cite{Feng} with $p = T$ and $n = N$. Note that the proof in \cite{Feng} also implies the bound (see equation S.82)

\begin{align} \label{N-negli}
    \limsup_{n \to \infty} N(n,k) \leq \frac{C}{k!}
\end{align}
where $C$ is a constant free of $n,p$ and $k$. $\hfill$ $\square$

\begin{lemma} \label{dterm}
Recall $\triangle_n$, $A_n$, $C_I$, $D(n,j)$ defined in (\ref{def-triang}), (\ref{def-An}), (\ref{def-CI}), and (\ref{def-Dnj}) respectively. Then, for every fixed $k \in \mb{N}$, we have 
$$\lim_{n \to \infty} D(n,k) = 0.$$
\end{lemma}

 {\noindent \textbf{Proof of Lemma \ref{dterm}}.}
    Fix $\ve>0$ and $k \in \mb{N}$. Define
\begin{align*}
    A^{+}_{n,\ve} &:= \la R_n \leq x+\ve, B_n \leq y+ \ve \ra, \\
    A^{-}_{n, \ve} &:= \la R_n \geq x-\ve, B_n \geq y- \ve \ra.
\end{align*}
Recall $N(n,k)$ as defined in (\ref{def-Nnk}). We have the bounds
\begin{align*}
    D(n,k) & \geq \sum_{I_1<I_2<\cdots<I_k} \Big[ \mb{P} \lb A^{-}_{n,\ve} C_{I_1} C_{I_2} \cdots C_{I_{k}} \rb - \mb{P} \lb A^{-}_{n,\ve} \rb \cdot \mb{P} \lb C_{I_1} C_{I_2} \cdots C_{I_{k}} \rb \Big] \\
    & - \left[ \mb{P} \lb A_n \rb - \mb{P} \lb A^{-}_{n,\ve} \rb \right] \cdot \sum_{I_1<I_2<\cdots<I_j}  \mb{P} \lb C_{I_1} C_{I_2} \cdots C_{I_{k}} \rb \\
    &= \sum_{I_1<I_2<\cdots<I_k} \Big[ \mb{P} \lb A^{-}_{n,\ve} C_{I_1} C_{I_2} \cdots C_{I_{k}} \rb - \mb{P} \lb A^{-}_{n,\ve} \rb \cdot \mb{P} \lb C_{I_1} C_{I_2} \cdots C_{I_{k}} \rb \Big] \\
    & -  \left[ \mb{P} \lb A_n \rb - \mb{P} \lb A^{-}_{n,\ve} \rb \right] \cdot N(n,k),
\end{align*}
and similarly,
\begin{align*}
    D(n,k) & \leq \sum_{I_1<I_2<\cdots<I_k} \mb{P} \lb A^{+}_{n,\ve} C_{I_1} C_{I_2} \cdots C_{I_{k}} \rb - \mb{P} \lb A^{+}_{n,\ve} \rb \cdot \mb{P} \lb C_{I_1} C_{I_2} \cdots C_{I_{k}} \rb \\
    & + \left[  \mb{P} \lb A^{+}_{N,\ve}  \rb - \mb{P} \lb A_n \rb  \right] \cdot N(n,k).
\end{align*}
Due to the estimate (\ref{N-negli}) below, and by letting $\ve \to 0$, it is sufficient to show that 
\begin{align}
  \lim_{\ve \to 0} \liminf_{n \to \infty}  \sum_{I_1<I_2<\cdots<I_k} \Big[ \mb{P} \lb A^{-}_{n,\ve} C_{I_1} C_{I_2} \cdots C_{I_{k}} \rb - \mb{P} \lb A^{-}_{n,\ve} \rb \cdot \mb{P} \lb C_{I_1} C_{I_2} \cdots C_{I_{k}} \rb \Big]  \to 0 \label{upperD}
\end{align}
and 
\begin{align}
   \lim_{\ve \to 0} \limsup_{n \to \infty}  \sum_{I_1<I_2<\cdots<I_k} \Big[ \mb{P} \lb A^{+}_{n,\ve} C_{I_1} C_{I_2} \cdots C_{I_{k}} \rb - \mb{P} \lb A^{+}_{n,\ve} \rb \cdot \mb{P} \lb C_{I_1} C_{I_2} \cdots C_{I_{k}} \rb \Big]  \to 0 \label{lowerD}.
\end{align}
Since (\ref{upperD}) and (\ref{lowerD}) are similar, we will only show (\ref{lowerD}). Fix a configuration $I_1<I_2<\cdots<I_k$ and let $I_l=(a_l,b_l)$, $l=1,2,\cdots,k$. Define
\begin{align}
    \triangle_{n,k} &:= \la (i_l,j): i_l +1 \leq j \leq n, 1 \leq l \leq k \ra \cup \la (i,j_l): 1 \leq i \leq j_l-1, 1 \leq l \leq k \ra, \label{def-triangNK}\\
    R^{'}_n &:= \frac{\sqrt{2p}}{n} \sum_{I=(i,j) \in \triangle_{n,k}}  \bm{X}^{\top}_i \bm{X}_j, \label{def-R'} \\
    B^{'}_n &:= \frac{p}{n} \sum_{I=(i,j) \in \triangle_{n,k}} \left[  \lb \bm{X}^{\top}_i \bm{X}_j \rb^2 -\frac{1}{p} \right]. \label{def-B'}
\end{align}
Recall $R_n$ and $B_n$ being the Rayleigh test and Bingham test defined in \eqref{Rayleigh} and \eqref{Bingham}, respectively. It is easy to see that $R_n-R^{'}_n$ and $B_n-B^{'}_n$ are independent from the events $\la C_{I_l} \ra_{l=1}^{k}$. Write

\begin{align*}
    \mb{P} \lb A^{+}_{n,\ve} C_{I_1} C_{I_2} \cdots C_{I_{k}} \rb &= \mb{P} \lb R_n \leq x+\ve, B_n \leq y+ \ve ,C_{I_1} C_{I_2} \cdots C_{I_{k}} \rb \\
    &=  \mb{P} \lb R_n \leq x+\ve, |R^{'}_n| \leq \ve, B_n \leq y+ \ve , |B^{'}_n| \leq \ve, C_{I_1} C_{I_2} \cdots C_{I_{k}} \rb \\
    &+ \mb{P} \lb \max \la |R^{'}_n| , |B^{'}_n|   \ra >\ve \rb \\
    & \leq   \mb{P} \lb R_n -R^{'}_n \leq x+2\ve, B_n -B^{'}_n \leq y+ 2\ve, C_{I_1} C_{I_2} \cdots C_{I_{k}} \rb \\
    & + \mb{P} \lb \max \la |R^{'}_n| , |B^{'}_n|   \ra >\ve \rb \\
    & = \mb{P} \lb R_n -R^{'}_n \leq x+2\ve, B_n -B^{'}_n \leq y+ 2\ve \rb \cdot \mb{P} \lb C_{I_1} C_{I_2} \cdots C_{I_{k}} \rb \\
    & + \frac{2C(\tau,k)}{\ve^{2\tau}n^\tau},
\end{align*}
where we used Lemma \ref{swapping} in the last inequality, and $C(\tau,k)$ is the constant in this lemma. Moreover, Lemma \ref{swapping} again gives 
\begin{align*}
   & \mb{P} \lb R_n -R^{'}_n \leq x+2\ve, B_n -B^{'}_n \leq y+ 2\ve \rb  \\
   \leq  & \mb{P} \lb R_n -R^{'}_n \leq x+2\ve, |R^{'}_n| \leq \ve, B_n -B^{'}_n \leq y+ 2\ve, |B^{'}_n| \leq \ve \rb 
   + \mb{P} \lb \max \la |R^{'}_n| , |B^{'}_n|   \ra >\ve \rb  \\
   \leq &  \mb{P} \lb A^{+}_{n,3\ve} \rb + \frac{2C(\tau,k)}{\ve^{2\tau} n^\tau}.
\end{align*}
Consequently,
\begin{align*}
        \mb{P} \lb A^{+}_{n,\ve} C_{I_1} C_{I_2} \cdots C_{I_{k}} \rb \leq \left[   \mb{P} \lb A^{+}_{n,3\ve} \rb + \frac{2C(\tau,k)}{\ve^{2\tau} n^\tau} \right] \cdot \mb{P} \lb C_{I_1} C_{I_2} \cdots C_{I_{k}} \rb + \frac{2C(\tau,k)}{\ve^{2\tau} n^\tau}.
\end{align*}
Sum up over all configurations $I_1<I_2<\cdots<I_k$ in $\triangle_n$, we have
\begin{align*} 
    & \sum_{I_1<I_2<\cdots<I_k} \Big[ \mb{P} \lb A^{+}_{n,\ve} C_{I_1} C_{I_2} \cdots C_{I_{k}} \rb - \mb{P} \lb A^{+}_{n,\ve} \rb \cdot \mb{P} \lb C_{I_1} C_{I_2} \cdots C_{I_{k}} \rb \Big] \\
    \leq & \left[\mb{P} \lb A^{+}_{n,3\ve} \rb - \mb{P} \lb A^{+}_{n,\ve} \rb + \frac{2C(\tau,k)}{\ve^{2\tau} n^\tau} \right] \cdot N(n,k) + \Big| \triangle_n \Big| \cdot \frac{2C(\tau,k)}{\ve^{2\tau} n^\tau}.
\end{align*}
Take $\tau=2k+1$ and use the bounds (\ref{N-negli}) together with $|\triangle_n| \leq n^{2k}$, gives
\begin{align*}
 & \limsup_{n \to \infty}  \sum_{I_1<I_2<\cdots<I_k} \Big[ \mb{P} \lb A^{+}_{n,\ve} C_{I_1} C_{I_2} \cdots C_{I_{k}} \rb - \mb{P} \lb A^{+}_{n,\ve} \rb \cdot \mb{P} \lb C_{I_1} C_{I_2} \cdots C_{I_{k}} \rb \Big] \\
 \leq & \limsup_{n \to \infty}\frac{C \left[ \mb{P} \lb A^{+}_{n,3\ve} \rb - \mb{P} \lb A^{+}_{n,\ve} \rb \right]}{k!},
\end{align*}
where $C$ is the constant in Lemma \ref{Nterm}. Note that we have also used Lemma \ref{Nterm} to get the bound of order $1/k!$ on $N(n,k)$ in the above. By sending $\ve \to 0$, we get (\ref{lowerD}) since both $\mb{P} \lb A^{+}_{n,3\ve} \rb$ and $\mb{P} \lb A^{+}_{n,\ve} \rb$ converge to $\Phi(x)\Phi(y)$. This concludes the proof. $\hfill$ $\square$

\begin{lemma} \label{swapping}
    Recall $\triangle_{nk}$, $R^{'}_n$ and $B^{'}_n$ defined in (\ref{def-triangNK}), (\ref{def-R'}) and (\ref{def-B'}). For any $\ve>0$ and any integer $\tau>0$, we have the concentration inequalities
    \begin{align*}
     \mb{P} \lb |  R^{'}_{n} | \geq \ve \rb & \leq \frac{C(\tau,k)}{\ve^{2 \tau} n^{\tau}}, \\
     \mb{P} \lb | B^{'}_{n} | \geq \ve \rb & \leq \frac{C(\tau,k)}{\ve^{2\tau} n^{\tau}},
    \end{align*}
    where $C(\tau,k)$ is a constant depending only on $\tau$ and $k$. 
\end{lemma}

 {\noindent \textbf{Proof of Lemma \ref{swapping}}.}
    The two inequalities are simple applications of Markov inequality and a moment bound. Let us provide the details below. It is easy to see that the size of the set $\triangle_{n,k}$ is at most 
$$\sum_{l=1}^{k} n - i_l + 1 + n-j_l+1 \leq 2nk + 2 \leq 3nk.$$
We claim that
\begin{align}
    \mb{E} \left[ \sqrt{p} \sum_{I=(i,j) \in \triangle_{n,k}} \bm{X}^{\top}_i \bm{X}_j \right]^{2 \tau} &\leq C(\tau,k) \cdot n^{\tau} ; \label{moment-bound1} \\
      p^{2\tau} \cdot  \mb{E} \left[  \sum_{I=(i,j) \in \triangle_{n,k}} \lb \bm{X}^{\top}_i \bm{X}_j \rb^2 -\frac{1}{p} \right]^{2 \tau} & \leq C(\tau,k) \cdot n^{\tau} \label{moment-bound2},
\end{align}
for some constant $C(\tau,k)$ depending only on $k$ and $\tau$. We show them next.

To show \eqref{moment-bound1} and \eqref{moment-bound2}, consider the following grouping rule, which partition $\triangle_{n,k}$ into disjoint union of  $2k$ groups, say $G_1, G_2, \cdots, G_{2k}$. Starting with $i_1$, put all the pairs $\la (i_1,j): I=(i_1,j) \in \triangle_{n,k} \ra$ into $G_1$. Next, consider $\triangle_{n,k} \setminus G_1$ and repeat the same procedure with index $i_2$. Once we finish with all indexs $i_1, i_2, \cdots, i_k$, the $k+1$-th group will start with index $j_1$ and continue until the last group, $G_{2k}$, is formed. This create a partition of $\triangle_{n,k}$ since all pairs $(i,j)$ belong to $\triangle_{n,k}$ either has $i=i_l$ for some $1 \leq l \leq k$ or $j=j_l$ for some $1 \leq l \leq k$.

It may happen that a few of these groups are empty, but that does not affect the argument below. The only two things we need are 
\begin{align*}
    \triangle_{n,k} =  \bigcup_{j=1}^{2k} G_j \quad \text{and} \quad G_i \cap G_j = \emptyset
\end{align*}
for all $i \neq j$. Thanks to the elementary inequality $\lb \sum_{i=1}^{n} a_i \rb^{\tau} \leq n^{\tau-1} \sum_{i=1}^{n} |a_i|^\tau$, we get

\begin{align*}
        \mb{E} \left[ \sqrt{p} \sum_{I=(i,j) \in \triangle_{n,k}} \bm{X}^{\top}_i \bm{X}_j \right]^{2 \tau} &\leq (2k)^{2\tau-1} \cdot \sum_{l=1}^{2k} \mb{E} \left[ \sqrt{p} \sum_{I=(i,j) \in G_l} \bm{X}^{\top}_i \bm{X}_j \right]^{2 \tau} \\
        &\leq (2k)^{2\tau-1} \cdot  \sum_{l=1}^{2k} |G_l|^{\tau} \cdot p^{\tau} \cdot \mb{E} \lb \bm{X}^{\top}_1 \bm{X}_2 \rb^{2 \tau} \\
        & \leq (2k)^{2\tau} \cdot \lb \sum_{l=1}^{2k} |G_l| \rb^{\tau} \cdot p^{\tau} \cdot \mb{E} \lb \bm{X}^{\top}_1 \bm{X}_2 \rb^{2 \tau},
\end{align*}
where we used the moment inequality for sum of independent random variables (conditionally on the common index in each group). Since $|G_1|+|G_2|+\cdots+|G_{2k}| = |\triangle_{nk}|\leq 3nk$, the last bound simplies to

\begin{align*}
     \mb{E} \left[ \sqrt{p} \sum_{I=(i,j) \in \triangle_{n,k}} \bm{X}^{\top}_i \bm{X}_j \right]^{2 \tau} &\leq 3^{\tau} \cdot (2k)^{2\tau} \cdot k^{\tau} \cdot n^{\tau} \cdot \sup_{p \geq 1} \left[ p^{\tau} \cdot \mb{E} \lb \bm{X}^{\top}_1 \bm{X}_2 \rb^{2 \tau} \right].
\end{align*}
It suffices to check that the last supremum is finite and depends only on $\tau$. To see this, by Lemma 2.4 in \cite{Jiang09}, we have
\begin{align*}
    \sup_{p \geq 1} \left[ p^{\tau} \cdot \mb{E} \lb \bm{X}^{\top}_1 \bm{X}_2 \rb^{2 \tau} \right] &= \sup_{p \geq 1} \left[ p^{\tau} \cdot \frac{(2\tau-1)!!}{\Pi_{i=1}^{\tau} (p+2i-2)} \right] = C(\tau)<\infty,
\end{align*}
which depends only on $r$. This completes the proof of (\ref{moment-bound1}). The proof of (\ref{moment-bound2}) is the same with the only difference being that we have to check the finiteness of 
\begin{align*}
       & \sup_{p \geq 1} \left[ p^{2\tau} \cdot \mb{E} \lb \lb \bm{X}^{\top}_1 \bm{X}_2 \rb^2 -\frac{1}{p} \rb^{2\tau} \right] \\
        \leq &  2^{2\tau-1} \cdot \sup_{p \geq 1} \left[ p^{2\tau} \cdot  \mb{E} \lb \bm{X}^{\top}_1 \bm{X}_2 \rb^{4\tau} + 1 \right] \\
        = & 2^{2\tau-1} \cdot \sup_{p \geq 1} \left[ p^{2\tau} \cdot \frac{(2\tau-1)!!}{\Pi_{i=1}^{2\tau} (p+2i-2)} + 1\right],
\end{align*}
which is indeed finite and depend only on $\tau$. Note that use  Lemma 2.4 in \cite{Jiang09} to get the value of the expectation in the above.  Finally, by using the Markov inequality and the bound (\ref{moment-bound1}), we obtain 
\begin{align*}
    \mb{P} \lb |  R^{'}_{n} | \geq \ve \rb &\leq \frac{\mb{E} \lb |R^{'}_n|^{2 \tau} \rb}{\ve^{2 \tau}} \\
    & \leq \frac{2^\tau }{n^{2\tau} \ve^{2\tau}}  \cdot p^{\tau}  \mb{E} \Big( \sum_{I=(i,j) \in \triangle_{n,k}} \bm{X}^{\top}_i \bm{X}_j \Big)^{2\tau} \\
    & \leq \frac{2^\tau}{n^{2\tau} \ve^{2\tau}}  \cdot C(\tau,k) \cdot n^{\tau} = \frac{C^{'}(\tau,k)}{n^\tau \ve^{2\tau}}.
\end{align*}
The concentration inequality for $B_n$ can be proven similarly by using the Markov inequality and the moment bound (\ref{moment-bound2}). The proof is completed. $\hfill$ $\square$

\bibliographystyle{plain}
\bibliography{IAI}

\begin{thebibliography}{10}

\bibitem{albisetti2020testing}
I.~Albisetti, F.~Balabdaoui, and H.~Holzmann.
\newblock Testing for spherical and elliptical symmetry.
\newblock {\em Journal of Multivariate Analysis}, 180:104667, 2020.

\bibitem{complex}
J.~Bak and D.J. Newman.
\newblock {\em Complex Analysis}.
\newblock Springer New York, 2010.

\bibitem{Banerjee03}
A.~Banerjee, I.~Dhillon, J.~Ghosh, and S.~Sra.
\newblock Generative model-based clustering of directional data.
\newblock {\em Proceedings of the ninth ACM SIGKDD international conference on Knowledge discovery and data mining}, pages 19--28, 2003.

\bibitem{Banerjee04}
A.~Banerjee and J.~Ghosh.
\newblock Frequency sensitive competitive learning for scalable balanced clustering on high-dimensional hyperspheres.
\newblock {\em IEEE T. Neural Network.}, 15:702--719, 2004.

\bibitem{Jiang13}
T.~Cai, J.~Fan, and T.~Jiang.
\newblock Distributions of angles in random packing on spheres.
\newblock {\em J. Mach. Learn. Res.}, 14:1837--1864, 2013.

\bibitem{CH09}
Y-B. Chan and P.~Hall.
\newblock Robust nearest-neighbor methods for classifying high-dimensional data.
\newblock {\em Ann. Statist. 37(6A): 3186-3203.}, 37(6A):3186--3203, 2009.

\bibitem{Cutting-P-V}
C.~Cutting, D.~Paindaveine, and T.~Verdebout.
\newblock Testing uniformity on high-dimensional spheres against contiguous rotationally symmetric alternatives.
\newblock {\em Ann. Stat.}, 45:1024--1058, 2017.

\bibitem{Cutting-P-V2}
C.~Cutting, D.~Paindaveine, and T.~Verdebout.
\newblock Testing uniformity on high-dimensional spheres: The non-null behaviour of the bingham test.
\newblock {\em Annales de I'I.H.P. Probabilit\'es et statistiques}, 58:567--602, 2022.

\bibitem{Darling}
D.~A. Darling.
\newblock The influence of the maximum term in the addition of independent random variables.
\newblock {\em Transactions of the American Mathematical Society}, 73(1):95--107, 1952.

\bibitem{Dryden05}
I.~L. Dryden.
\newblock Statistical analysis on high-dimensional spheres and shape spaces.
\newblock {\em Ann. Stat.}, 33:1643--1665, 2005.

\bibitem{Feng}
L.~Feng, T.~Jiang, B.~Liu, and W.~Xiong.
\newblock Max-sum tests for cross-sectional independence of high-demensional panel data.
\newblock {\em Ann. Stat. 50(2), 1124-1143.}, 2022.

\bibitem{fernandez2023new}
A.~Fern{\'a}ndez-de Marcos and E.~Garc{\'\i}a-Portugu{\'e}s.
\newblock On new omnibus tests of uniformity on the hypersphere.
\newblock {\em Test}, 32(4):1508--1529, 2023.

\bibitem{Fisher}
N.~I. Fisher, T.~Lewis, and B.~J. Embleton.
\newblock {\em Statistical Analysis of Spherical Data}.
\newblock Cambridge Univ. Press press, 1987.

\bibitem{garcia2021cramer}
E.~Garc{\'\i}a-Portugu{\'e}s, P.~Navarro-Esteban, and J.~Cuesta-Albertos.
\newblock A cram{\'e}r--von mises test of uniformity on the hypersphere.
\newblock In {\em Statistical Learning and Modeling in Data Analysis: Methods and Applications 12}, pages 107--116. Springer, 2021.

\bibitem{garcia2023projection}
E.~Garc{\'\i}a-Portugu{\'e}s, P.~Navarro-Esteban, and J.~Cuesta-Albertos.
\newblock On a projection-based class of uniformity tests on the hypersphere.
\newblock {\em Bernoulli}, 29(1):181--204, 2023.

\bibitem{Hall}
P.~Hall and C.~C. Heyde.
\newblock {\em Martingale Limit Theory and Its Application}.
\newblock Academic Press, New York., 1980.

\bibitem{heiny2017eigenvalues}
J.~Heiny and T.~Mikosch.
\newblock Eigenvalues and eigenvectors of heavy-tailed sample covariance matrices with general growth rates: the iid case.
\newblock {\em Stochastic Processes and their Applications}, 127(7):2179--2207, 2017.

\bibitem{heiny2019eigenstructure}
J.~Heiny and T.~Mikosch.
\newblock The eigenstructure of the sample covariance matrices of high-dimensional stochastic volatility models with heavy tails.
\newblock {\em Bernoulli}, 25(4B):3590--3622, 2019.

\bibitem{heiny2021large}
J.~Heiny and T.~Mikosch.
\newblock Large sample autocovariance matrices of linear processes with heavy tails.
\newblock {\em Stochastic Processes and their Applications}, 141:344--375, 2021.

\bibitem{heiny2022}
J.~Heiny and J.~Yao.
\newblock Limiting distributions for eigenvalues of sample correlation matrices from heavy-tailed populations.
\newblock {\em Ann. Stat.}, 50(6):3249--3280, 2022.

\bibitem{Jiang09}
T.~Jiang.
\newblock A variance formula related to quantum conductance.
\newblock {\em Physics Letters A 373, 2117-2121.}, 2009.

\bibitem{jiang2023asymptotic}
T.~Jiang and T.~Pham.
\newblock Asymptotic distributions of largest pearson correlation coefficients under dependent structures.
\newblock {\em arXiv preprint arXiv:2304.13102}, 2023.

\bibitem{Juan2001}
J.~Juan and F.~J. Prieto.
\newblock Using angles to identify concentrated multivariate outliers.
\newblock {\em Technometrics}, 43:311--322, 2001.

\bibitem{Kulik2020}
R.~Kulik and P.~Soulier.
\newblock {\em Heavy-Tailed Time Series}.
\newblock Springer New York, 2020.

\bibitem{Ley-Verdebout}
C.~Ley and T.~Verdebout.
\newblock {\em Modern Directional Statistics}.
\newblock CRC Press, Boca Raton., 2017.

\bibitem{M-Jupp}
K.~V. Mardia and P.~E. Jupp.
\newblock {\em Directional Statistics}.
\newblock John Wiley \& Sons, Chichester., 2000.

\bibitem{Ley-P}
D.~Paindaveine and T.~Verdebout.
\newblock On high-dimensional sign tests.
\newblock {\em Bernoulli}, 22:1745--1769, 2016.

\bibitem{Ley-P-2}
D.~Paindaveine and T.~Verdebout.
\newblock Detecting the direction of a signal on high-dimensional spheres: non-null and le cam optimality results.
\newblock {\em Probability Theory and Related Fields 176, 1165-1216.}, 2020.

\bibitem{pewsey2021recent}
A.~Pewsey and E.~Garc{\'\i}a-Portugu{\'e}s.
\newblock Recent advances in directional statistics.
\newblock {\em Test}, 30(1):1--58, 2021.

\bibitem{survey-uni}
E.~G. Portugu\'es and T.~Verdebout.
\newblock An overview of uniformity tests on the hypersphere.
\newblock {\em arXiv:1804.00286.}, 2018.

\bibitem{shao1997self}
Qi-Man Shao.
\newblock Self-normalized large deviations.
\newblock {\em The Annals of Probability}, 25(1):285--328, 1997.

\bibitem{tang2023nonparametric}
Y.~Tang and B.~Li.
\newblock A nonparametric test for elliptical distribution based on kernel embedding of probabilities.
\newblock {\em arXiv preprint arXiv:2306.10594}, 2023.

\bibitem{WPL15}
L.~Wang, B.~Peng, and R.~Li.
\newblock A high-dimensional nonparametric multivariate test for mean vector.
\newblock {\em Journal of the American Statistical Association}, 110:1658--1669, 2015.

\bibitem{Yu1}
X.~Yu, D.~Li, and L.~Xue.
\newblock Fisher's combined probability test for high-dimensional covariance matrices.
\newblock {\em Journal of the American Statistical Association}, 119:511--524, 2022.

\bibitem{Yu2}
X.~Yu, D.~Li, L.~Xue, and R.~Li.
\newblock Power-enhanced simultaneous test of high-dimensional mean vectors and covariance matrices with application to gene-set testing.
\newblock {\em Journal of the American Statistical Association}, 118:2548--2561, 2023.

\bibitem{Yu3}
X.~Yu, J.~Yao, and L.~Xue.
\newblock Power enhancement for testing multi-factor asset pricing models via fisher's method.
\newblock {\em Journal of Econometrics}, 239:105458, 2024.

\bibitem{Zou14}
C.~Zou, L.~Peng, L.~Feng, and Z.~Wang.
\newblock Multivariate-sign-based high-dimensional tests for sphericity.
\newblock {\em Biometrika}, 101:229--236, 2014.

\end{thebibliography}

\end{document}